\documentclass[a4paper,12pt]{article}
\usepackage{amsfonts}
\usepackage{amsbsy}

\newcommand{\card}[1]{\left|#1\right|}
\newcommand{\setof}[1]{\left\{#1\right\}}
\newcommand{\R}{\mathbb{R}}
\newcommand{\Z}{\mathbb{Z}}
\newcommand{\Tr}{\mathop{\mathrm{Tr}}}

\newtheorem{thm}{Theorem}

\newtheorem{cor}{Corollary}
\newtheorem{cutset}{Cutset Lemma}
\newtheorem{proof}{Proof}

\newenvironment{pf}{\proof\normalfont}{\endproof}
\newcommand{\halmos}{\protect\nolinebreak\mbox{\quad\rule{1ex}{1ex}}}
\newtheorem{example}{Example}
\newenvironment{eg}{\example\normalfont}{\endexample}

\newcommand{\blobb}{\circle*{0.4}}
\newcommand{\blobbb}{\circle*{7}}
\newcommand{\xlabel}[2]{\put(#1,2.9){\line(0,1){0.2}}
  \put(#1,2.6){\makebox(0,0.3){$#2$}}}
\newcommand{\ylabel}[2]{\put(0.9,#1){\line(1,0){0.2}}
   \put(0.4,#1){\makebox(0.4,0)[r]{$#2$}}}
\newcommand{\xxxlabel}[2]{\put(#1,1.95){\line(0,1){0.1}}
  \put(#1,1.7){\makebox(0,0.3){$#2$}}}
\newcommand{\yyylabel}[2]{\put(-0.05,#1){\line(1,0){0.1}}
   \put(-0.5,#1){\makebox(0.4,0)[r]{$#2$}}}

\title{Using graphs to find the best block designs}

\author{R. A. Bailey and Peter J. Cameron}

\date{}

\begin{document}

\maketitle

R. A. Bailey obtained a DPhil in group theory from the University of 
Oxford. She worked at the Open University, and then held a
post-doctoral research fellowship in Statistics at the University of 
Edinburgh. This was followed by ten years in the Statistics Department
at Rothamsted Experimental Station, which at that time came under the
auspices of the
Agriculture and Food Research Council.  She returned to university
life as Professor of Mathematical Sciences at Goldsmiths' College, 
University of London, and has been Professor of Statistics at Queen
Mary, University of London, since 1994.

Peter J. Cameron is Professor of Mathematics at Queen Mary, University of
London, where he has been since 1986, following a position as tutorial
fellow at Merton College, Oxford. Since his DPhil in Oxford, he has been
interested in a variety of topics in algebra and combinatorics,
especially their interactions. He has held visiting positions at the
University of Michigan, California Institute of Technology, and the
University of Sydney. He is currently chair of the British Combinatorial
Committee.

\tableofcontents

\begin{abstract}
A statistician designing an experiment wants to get as much information as
possible from the data gathered. Often this means the most precise estimate
possible (that is, an estimate with minimum possible variance) of the
unknown parameters. If there are several parameters, this can be interpreted
in many ways: do we want to minimize the average variance, or the maximum
variance, or the volume of a confidence region for the parameters?

In the case of block designs, these optimality criteria can be calculated
from the concurrence graph of the design, and in many cases from its
Laplacian eigenvalues.  The Levi graph can also be used.
The various criteria turn out to be closely connected
with other properties of the graph as a network, such as number of spanning
trees, isoperimetric number, and the sum of the resistances between pairs
of vertices when the graph
is regarded as an electrical network.

In this chapter, we discuss the notions of optimality for incomplete-block
designs, explain the graph-theoretic connections, and prove some old and new
results about optimality.
\end{abstract}

\section{What makes an incomplete-block design good for experiments?}
\label{sec:intro}

Experiments are designed in many ways: for example, Latin squares,
block designs, split-plot designs.  Combinatorialists, on the other
hand, have a much more specialized usage of the term ``design'', as we
remark later.  We are concerned here with incomplete-block designs,
more special than the statistician's designs and more general than the
mathematician's. 

To a statistician, a \emph{block design} has two components.
There is an underlying set of experimental units, partitioned into
$b$~blocks of size~$k$.  There is a further set of $v$~treatments, and
also a function~$f$ from units to treatments, specifying which treatment 
is allocated to which experimental unit; that is, $f(\omega)$ is the 
treatment allocated to experimental unit~$\omega$.  Thus each block 
defines a subset, or maybe a multi-subset, of the treatments.

In a \emph{complete-block design}, we have $k=v$ and each treatment occurs 
once in every block.  Here we assume that blocks are \emph{incomplete} in the 
sense that $k<v$.

We assume that the purpose of the experiment is to find out about the 
treatments, and differences between them.  The blocks are an unavoidable 
nuisance, an inherent feature of the experimental units.  In an agricultural
experiment the experimental units may be field plots and the blocks may be 
fields or plough-lines;  in a clinical trial the experimental units may be
patients and the blocks hospitals; in process engineering the experimental 
units may be runs of a machine that is recalibrated each day and the 
blocks days.  See \cite{rabbook} for further examples.

In all of these situations, the values of $b$, $k$ and~$v$ are given.
Given these values, not all incomplete-block designs are 
equally good. This chapter describes some criteria that can be used to choose 
between them.

For example, Fig.~\ref{fig:queen} shows two block designs with
$v=15$, $b=7$ and $k=3$.  We use the convention that the treatments are 
labelled $1$, \ldots,~$v$, that columns represent blocks, and that the order 
of the entries in each column is not significant.  Where necessary, we use 
the notation $\Gamma_j$ to refer to the block which is shown as the 
$j$th column, for $j=1$, \ldots, $b$.

\begin{figure}
\begin{center}
\begin{tabular}{c@{\qquad\qquad}c}
$\begin{array}{|c|c|c|c|c|c|c|}
\hline
1 & 1 & 2 & 3 & 4 & 5 & 6\\
2 & 4 & 5 & 6 & 10 & 11 & 12\\
3 & 7 & 8 & 9 & 13 & 14 & 15\\
\hline
\end{array}
$
&
$\begin{array}{|c|c|c|c|c|c|c|}
\hline
1 & 1 & 1 & 1 & 1 & 1 & 1\\
2 & 4 & 6 & 8 & 10 & 12 & 14\\
3 & 5 & 7 & 9 & 11 & 13 & 15\\
\hline
\end{array}
$\\ \\
(a) & (b)
\end{tabular}
\end{center}
\caption{Two block designs with $v=15$, $b=7$ and $k=3$}
\label{fig:queen}
\end{figure}

The \emph{replication} $r_i$ of treatment $i$ is defined to be
$\card{f^{-1}(i)}$, which is the number of experimental units to which 
it is allocated.  For the design in Fig.~\ref{fig:queen}(a), 
$r_i \in \setof{1,2}$ for all~$i$.  As we see later, statisticians
tend to prefer designs in which all the replications are as equal as 
possible.  
If $r_i=r_j$ for $1\leq i<j\leq v$ then the design is 
\emph{equireplicate}: then the common value of $r_i$ is usually written 
as~$r$, and $vr=bk$.

The design in Fig.~\ref{fig:queen}(b) is a \emph{queen-bee design} 
because there is (at least) one treatment that occurs in every block.  
Scientists tend to prefer such designs because they have been taught to 
compare every treatment to one distinguished treatment, which may be called
a \emph{control treatment}.

\begin{figure}[htbp]
\begin{center}
\begin{tabular}{c@{\qquad\qquad}c}
$\begin{array}{|c|c|c|c|c|c|c|}
\hline
1 & 1 & 1 & 1 & 2 & 2 & 2\\
2 & 3 & 3 & 4 & 3 & 3 & 4\\
3 & 4 & 5 & 5 & 4 & 5 & 5\\
\hline
\end{array}
$
&
$\begin{array}{|c|c|c|c|c|c|c|}
\hline
1 & 1 & 1 & 1 & 2 & 2 & 2\\
1 & 3 & 3 & 4 & 3 & 3 & 4\\
2 & 4 & 5 & 5 & 4 & 5 & 5\\
\hline
\end{array}
$\\ \\
(a) & (b)
\end{tabular}
\end{center}
\caption{Two block designs with $v=5$, $b=7$ and $k=3$}
\label{fig:binary}
\end{figure}

Fig.~\ref{fig:binary} shows two block designs with $v=5$, $b=7$ and $k=3$.
The design in Fig.~\ref{fig:binary}(b) shows a new feature: 
treatment~$1$ occurs on two experimental units in block~$\Gamma_1$.  A 
block design is \emph{binary} if $f(\alpha)\ne f(\omega)$ whenever 
$\alpha$ and $\omega$ are experimental units in the same block.  The design in
Fig.~\ref{fig:binary}(a) is binary.  It seems to be obvious that 
binary designs must be better than non-binary ones, but we shall see later that
this is not necessarily so.  However, 
if there is any block on which $f$ is constant, 
then that block provides no information about
treatments, so we assume from now on that there are no such blocks.

\begin{figure}[htbp]
\begin{center}
\begin{tabular}{c@{\qquad\qquad}c}
$\begin{array}{|c|c|c|c|c|c|c|}
\hline
1 & 2 & 3 & 4 & 5 & 6 & 7\\
2 & 3 & 4 & 5 & 6 & 7 & 1\\
4 & 5 & 6 & 7 & 1 & 2 & 3\\
\hline
\end{array}
$
&
$\begin{array}{|c|c|c|c|c|c|c|}
\hline
1 & 2 & 3 & 4 & 5 & 6 & 7\\
2 & 3 & 4 & 5 & 6 & 7 & 1\\
3 & 4 & 5 & 6 & 7 & 1 & 2\\
\hline
\end{array}
$
\\ \\
(a) & (b)
\end{tabular}
\end{center}
\caption{Two block designs with $v=7$, $b=7$ and $k=3$}
\label{fig:fano}
\end{figure}

Fig.~\ref{fig:fano} shows two equireplicate binary block designs with 
$v=7$, $b=7$ and $k=3$.  A binary design is \emph{balanced} if every pair
of distinct treatments occurs together in the same number of blocks.
If that number is $\lambda$, then $r(k-1) = (v-1)\lambda$.  Such designs 
are also called \emph{$2$-designs} or \emph{BIBDs}.  The design in 
Fig.~\ref{fig:fano}(a) is balanced with $\lambda=1$; the design in 
Fig.~\ref{fig:fano}(b) is not balanced.  

Pure mathematicians usually assume that, if they exist,  balanced designs 
are better than non-balanced ones.  (Indeed, many do not call a structure a
`design' unless it is balanced.)
As we shall show in 
Section~\ref{sec:bibd}, this assumption is correct for all the 
criteria considered here.
However, for given values of $v$ and $k$, a non-balanced design with a larger
value of $b$ may produce more information than a balanced design with a 
smaller value of $b$.

\section{Graphs from block designs}

\subsection{The Levi graph}

\begin{figure}[hbtp]
\begin{center}
\setlength{\unitlength}{8mm}
\begin{picture}(12,3.5)(0,0.5)
\multiput(0,3)(2,0){7}{\blobb}
\multiput(2,1)(2,0){5}{\blobb}
\put(0,3.3){\makebox(0,0)[b]{$\Gamma_1$}}
\put(2,3.3){\makebox(0,0)[b]{$\Gamma_2$}}
\put(4,3.3){\makebox(0,0)[b]{$\Gamma_3$}}
\put(6,3.3){\makebox(0,0)[b]{$\Gamma_4$}}
\put(8,3.3){\makebox(0,0)[b]{$\Gamma_5$}}
\put(10,3.3){\makebox(0,0)[b]{$\Gamma_6$}}
\put(12,3.3){\makebox(0,0)[b]{$\Gamma_7$}}
\put(2,0.7){\makebox(0,0)[t]{1}}
\put(4,0.7){\makebox(0,0)[t]{2}}
\put(6,0.7){\makebox(0,0)[t]{3}}
\put(8,0.7){\makebox(0,0)[t]{4}}
\put(10,0.7){\makebox(0,0)[t]{5}}
\put(2,0.9){\line(-1,1){2}}
\put(2,1.1){\line(-1,1){2}}
\put(4,1){\line(-2,1){4}}
\put(2,1){\line(0,1){2}}
\put(6,1){\line(-2,1){4}}
\put(8,1){\line(-3,1){6}}
\put(2,1){\line(1,1){2}}
\put(6,1){\line(-1,1){2}}
\put(10,1){\line(-3,1){6}}
\put(2,1){\line(2,1){4}}
\put(8,1){\line(-1,1){2}}
\put(10,1){\line(-2,1){4}}
\put(4,1){\line(2,1){4}}
\put(6,1){\line(1,1){2}}
\put(8,1){\line(0,1){2}}
\put(4,1){\line(3,1){6}}
\put(6,1){\line(2,1){4}}
\put(10,1){\line(0,1){2}}
\put(4,1){\line(4,1){8}}
\put(8,1){\line(2,1){4}}
\put(10,1){\line(1,1){2}}
\end{picture} 
\end{center}
\caption{The Levi graph of the design in Fig.~\ref{fig:binary}(b)}
\label{fig:levibin}
\end{figure}
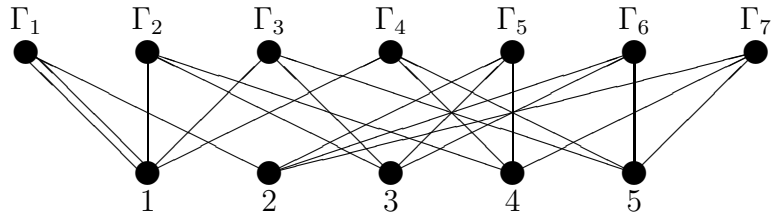

A simple way of representing a block design is its \emph{Levi graph}, 
or \textit{incidence graph},
introduced in \cite{levi}. This graph has $v+b$ vertices, one for each block 
and one for each treatment.  There are $bk$ edges, one for each experimental
unit.  If experimental unit $\omega$ is in block $j$ and $f(\omega) = i$, then
the corresponding edge~$\tilde e_\omega$ 
joins vertices $i$ and~$j$.  Thus the graph is 
bipartite, with one part consisting of block vertices and the other part 
consisting of treatment vertices. Moreover, the graph has multiple edges if
the design is not binary.
%
%
Fig.~\ref{fig:levibin} gives the Levi graph of the design in 
Fig.~\ref{fig:binary}(b).

We regard two block designs as the same if one can be obtained from the 
other by permuting 
the experimental units within each block.
Since the vertices of the Levi graph are labelled, there is a bijection
between block designs and their Levi graphs.

Let $n_{ij}$ be the number of edges from treatment-vertex $i$ to 
block-vertex~$j$; that is, treatment $i$ occurs on $n_{ij}$ experimental 
units in block $j$.  The $v\times b$ matrix $\mathbf{N}$ whose entries are the 
$n_{ij}$ is the \emph{incidence matrix} of the block design.  If the rows
and columns of $\mathbf{N}$ are labelled, we can recover the block design from
its incidence matrix.

\subsection{The concurrence graph}

In a binary design, the \emph{concurrence} $\lambda_{ij}$ of 
treatments~$i$ and $j$ is $r_i$ if $i=j$ and otherwise is the number 
of blocks in which $i$ and $j$ both occur.  For non-binary designs we have
to count the number of occurrences of the pair $\setof{i,j}$ in blocks 
according to multiplicity, so that $\lambda_{ij}$ is the $(i,j)$-entry of
$\boldsymbol{\Lambda}$, where $\boldsymbol{\Lambda} = \mathbf{N}
\mathbf{N}^\top$.   The matrix $\boldsymbol{\Lambda}$ is called  
the \emph{concurrence matrix} of the design.

The \emph{concurrence graph} of the design has the treatments as vertices.
There are no loops.  If $i\ne j$, then there are $\lambda_{ij}$ edges between
vertices $i$ and $j$.  Each such edge corresponds to a pair 
$\setof{\alpha,\omega}$ of experimental units in the same block, with
$f(\alpha)=i$ and $f(\omega)=j$: we denote this edge by 
$e_{\alpha\omega}$. (This edge does not join the experimental units $\alpha$
and $\omega$; it joins the treatments applied to these units.)
It follows that the degree $d_i$ of vertex $i$ is given by
\begin{equation}
d_i = \sum_{j\ne i} \lambda_{ij}.
\label{eq:valency}
\end{equation}
Figs.~\ref{fig:cgqueen} and~\ref{fig:cgbin} show the concurrence graphs of 
the designs in Figs.~\ref{fig:queen} and~\ref{fig:binary}, respectively.

\begin{figure}[htbp]
\begin{center}
\begin{tabular}{c@{\qquad\qquad\qquad}c}
\setlength{\unitlength}{8mm}
\begin{picture}(5,6)(0.8,0)
\put(0.4,4.6){\blobb}
\put(0.4,4.9){\makebox(0,0)[b]{13}}
\put(2.4,4.6){\blobb}
\put(2.4,4.9){\makebox(0,0)[b]{7}}
\put(3.6,4.6){\blobb}
\put(3.6,4.9){\makebox(0,0)[b]{8}}
\put(5.6,4.6){\blobb}
\put(5.6,4.9){\makebox(0,0)[b]{11}}
\put(0.4,3.4){\blobb}
\put(0.4,3.1){\makebox(0,0)[t]{10}}
\put(2.4,3.4){\blobb}
\put(2.1,3.1){\makebox(0,0)[rt]{1}}
\put(3.6,3.4){\blobb}
\put(3.9,3.1){\makebox(0,0)[tl]{2}}
\put(5.6,3.4){\blobb}
\put(5.6,3.1){\makebox(0,0)[t]{14}}
\put(1.2,1.4){\blobb}
\put(0.9,1.4){\makebox(0,0)[r]{12}}
\put(2.4,1.4){\blobb}
\put(2.6,1.2){\makebox(0,0)[tl]{6}}
\put(3.6,1.4){\blobb}
\put(3.9,1.4){\makebox(0,0)[l]{9}}
\put(1.8,0.4){\blobb}
\put(1.8,0.1){\makebox(0,0)[t]{15}}
\put(1.4,4.0){\blobb}
\put(1.4,4.3){\makebox(0,0)[b]{4}}
\put(4.6,4.0){\blobb}
\put(4.6,4.3){\makebox(0,0)[b]{5}}
\put(3.0,2.4){\blobb}
\put(3.3,2.4){\makebox(0,0)[l]{3}}
\put(1.2,1.4){\line(1,0){2.4}}
\put(1.8,0.4){\line(3,5){1.2}}
\put(1.8,0.4){\line(-3,5){0.6}}
\put(2.4,1.4){\line(3,5){1.2}}
\put(3.6,1.4){\line(-3,5){1.2}}
\put(2.4,3.4){\line(1,0){1.2}}
\put(2.4,3.4){\line(0,1){1.2}}
\put(2.4,3.4){\line(-5,3){2.0}}
\put(0.4,3.4){\line(0,1){1.2}}
\put(0.4,3.4){\line(5,3){2.0}}
\put(5.6,3.4){\line(0,1){1.2}}
\put(5.6,3.4){\line(-5,3){2.0}}
\put(3.6,3.4){\line(5,3){2}}
\put(3.6,3.4){\line(0,1){1.2}}
\end{picture}
&
\setlength{\unitlength}{8mm}
\begin{picture}(5,6)(0,0)
\put(2.5,2.5){\blobb}
\put(2.5,2.2){\makebox(0,0)[t]{1}}
\put(2.5,2.5){\line(5,1){2.5}}
\put(2.5,2.5){\line(3,2){2.1}}
\put(2.5,2.5){\line(2,3){1.4}}
\put(2.5,2.5){\line(1,5){0.5}}
\put(2.5,2.5){\line(-5,1){2.5}}
\put(2.5,2.5){\line(-3,2){2.1}}
\put(2.5,2.5){\line(-2,3){1.4}}
\put(2.5,2.5){\line(-1,5){0.5}}
\put(2.5,2.5){\line(-3,-2){2.1}}
\put(2.5,2.5){\line(-2,-3){1.4}}
\put(2.5,2.5){\line(-5,-1){2.5}}
\put(2.5,2.5){\line(5,-1){2.5}}
\put(2.5,2.5){\line(3,-2){2.1}}
\put(2.5,2.5){\line(2,-3){1.4}}
\put(5,3){\blobb}
\put(5.3,3){\makebox(0,0)[l]{12}}
\put(4.6,3.9){\blobb}
\put(4.9,3.9){\makebox(0,0)[l]{11}}
\put(3.9,4.6){\blobb}
\put(3.9,4.9){\makebox(0,0)[b]{10}}
\put(3,5){\blobb}
\put(3,5.3){\makebox(0,0)[b]{9}}
\put(0,3){\blobb}
\put(-0.3,3){\makebox(0,0)[r]{5}}
\put(0.4,3.9){\blobb}
\put(0.1,3.9){\makebox(0,0)[r]{6}}
\put(1.1,4.6){\blobb}
\put(1.1,4.9){\makebox(0,0)[b]{7}}
\put(2,5){\blobb}
\put(2,5.3){\makebox(0,0)[b]{8}}
\put(4.6,1.1){\blobb}
\put(4.9,1.1){\makebox(0,0)[l]{14}}
\put(3.9,0.4){\blobb}
\put(3.9,0.1){\makebox(0,0)[t]{15}}
\put(5,2){\blobb}
\put(5.3,2){\makebox(0,0)[l]{13}}
\put(0,2){\blobb}
\put(-0.3,2){\makebox(0,0)[r]{4}}
\put(0.4,1.1){\blobb}
\put(0.1,1.1){\makebox(0,0)[r]{3}}
\put(1.1,0.4){\blobb}
\put(1.1,0.1){\makebox(0,0)[t]{2}}
\put(3.9,0.4){\line(1,1){0.7}}
\put(5,2){\line(0,1){1}}
\put(4.6,3.9){\line(-1,1){0.7}}
\put(2,5){\line(1,0){1}}
\put(1.1,0.4){\line(-1,1){0.7}}
\put(0,2){\line(0,1){1}}
\put(0.4,3.9){\line(1,1){0.7}}
\end{picture}
\\ \\
(a) & (b)
\end{tabular}
\end{center}
\caption{The concurrence graphs of the designs in Fig.~\ref{fig:queen}}
\label{fig:cgqueen}
\end{figure}
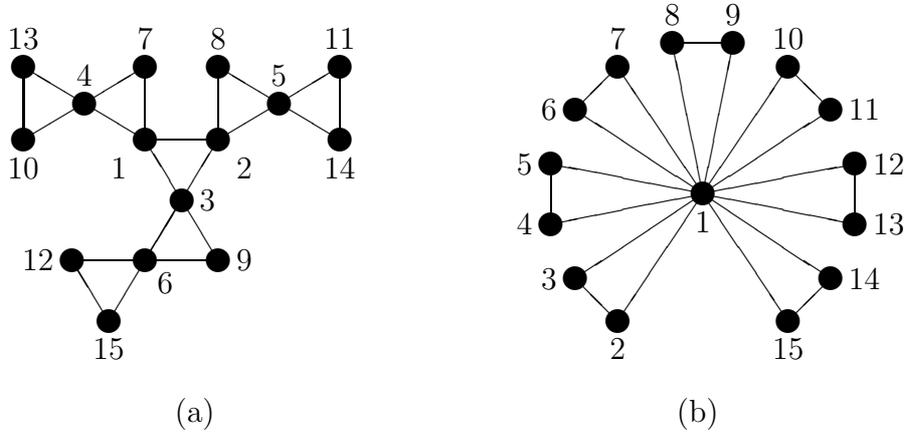

\begin{figure}
\begin{center}
\begin{tabular}{c@{\qquad\qquad\qquad}c}
\setlength{\unitlength}{8mm}
\begin{picture}(4,3.6)(-2,0)
\put(-1,0){\blobb}
\put(1,0){\blobb}
\put(-2,2){\blobb}
\put(2,2){\blobb}
\put(0,3){\blobb}
\put(-1,-0.1){\line(1,0){2}}
\put(-1,0.1){\line(1,0){2}}
\put(-2,2){\line(1,0){4}}
\put(-2,2){\line(2,1){2}}
\put(2,2){\line(-2,1){2}}
\put(-1.1,0.1){\line(1,3){1}}
\put(-0.9,-0.1){\line(1,3){1}}
\put(1.1,0.1){\line(-1,3){1}}
\put(0.9,-0.1){\line(-1,3){1}}
\put(-1.1,-0.1){\line(-1,2){1}}
\put(-0.9,0.1){\line(-1,2){1}}
\put(1.1,-0.1){\line(1,2){1}}
\put(0.9,0.1){\line(1,2){1}}
\put(-2.1,2.1){\line(2,1){2}}
\put(-1.9,1.9){\line(2,1){2}}
\put(2.1,2.1){\line(-2,1){2}}
\put(1.9,1.9){\line(-2,1){2}}
\put(-1.1,0.1){\line(3,2){3}}
\put(-0.9,-0.1){\line(3,2){3}}
\put(1.1,0.1){\line(-3,2){3}}
\put(0.9,-0.1){\line(-3,2){3}}
\put(0,3.3){\makebox(0,0)[b]{3}}
\put(2.3,2){\makebox(0,0)[l]{2}}
\put(1.3,0){\makebox(0,0)[l]{5}}
\put(-2.3,2){\makebox(0,0)[r]{1}}
\put(-1.3,0){\makebox(0,0)[r]{4}}
\end{picture}
&
\setlength{\unitlength}{8mm}
\begin{picture}(4,3)(-2,0)
\put(-1,0){\blobb}
\put(1,0){\blobb}
\put(-2,2){\blobb}
\put(2,2){\blobb}
\put(0,3){\blobb}
\put(-1,-0.1){\line(1,0){2}}
\put(-1,0.1){\line(1,0){2}}
\put(-2,2.1){\line(1,0){4}}
\put(-2,1.9){\line(1,0){4}}
\put(-1.1,0.1){\line(1,3){1}}
\put(-0.9,-0.1){\line(1,3){1}}
\put(1.1,0.1){\line(-1,3){1}}
\put(0.9,-0.1){\line(-1,3){1}}
\put(-1.1,-0.1){\line(-1,2){1}}
\put(-0.9,0.1){\line(-1,2){1}}
\put(1.1,-0.1){\line(1,2){1}}
\put(0.9,0.1){\line(1,2){1}}
\put(-2.1,2.1){\line(2,1){2}}
\put(-1.9,1.9){\line(2,1){2}}
\put(2.1,2.1){\line(-2,1){2}}
\put(1.9,1.9){\line(-2,1){2}}
\put(-1.1,0.1){\line(3,2){3}}
\put(-0.9,-0.1){\line(3,2){3}}
\put(1.1,0.1){\line(-3,2){3}}
\put(0.9,-0.1){\line(-3,2){3}}
\put(0,3.3){\makebox(0,0)[b]{3}}
\put(2.3,2){\makebox(0,0)[l]{2}}
\put(1.3,0){\makebox(0,0)[l]{5}}
\put(-2.3,2){\makebox(0,0)[r]{1}}
\put(-1.3,0){\makebox(0,0)[r]{4}}
\end{picture}
\\ \\
(a) & (b)
\end{tabular}
\end{center}
\caption{The concurrence graphs of the designs in Fig.~\ref{fig:binary}}
\label{fig:cgbin}
\end{figure}
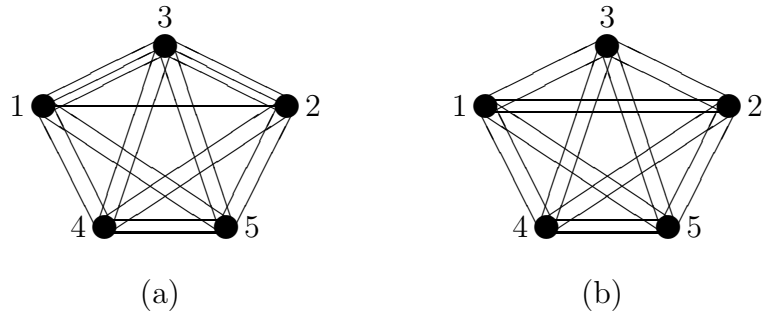

\begin{table}[hbtp]
\begin{center}
\begin{tabular}{c@{\quad\qquad}c}
$\left[
\begin{array}{rrrrr}
8 & -1 & -3 & -2 & -2\\
-1 & 8 & {-3} & -2 & -2\\
-3 & {-3} & 10 & -2 & -2\\
-2 & -2 & -2 & 8 & -2\\
-2 & -2 & -2 & -2 & 8
\end{array}
\right]$
&
$\left[
\begin{array}{rrrrr}
8 & {-2} & -2 & -2 & -2\\
{-2} & 8 & {-2} & -2 & -2\\
-2 & {-2} & 8 & -2 & -2\\
-2 & -2 & -2 & 8 & -2\\
-2 & -2 & -2 & -2 & 8
\end{array}
\right]
$
\\ \\
(a) & (b)
\end{tabular}
\end{center}
\caption{The Laplacian matrices of the concurrence graphs in 
Fig.~\ref{fig:cgbin}}
\label{tab:binary}
\end{table}

If $k=2$, then the concurrence graph is effectively the same as the block 
design. Although the block design cannot be recovered from the concurrence 
graph for larger values of $k$, we shall see in Section~\ref{sec:crit} that 
the concurrence graphs contain enough information to decide between two block
designs on any of the usual statistical criteria.  They were introduced
as \emph{variety concurrence graphs} in \cite{hdperw}, but 
are so useful that 
they may have been considered earlier.

\subsection{The Laplacian matrix of a graph}

Let $H$ be an arbitrary graph with $n$ vertices: it may have multiple
edges, but no loops. 
The 
\emph{Laplacian matrix} $\mathbf{L}$ of $H$ is defined to be the 
square matrix with 
rows and columns indexed by the vertices of $H$ whose $(i,i)$-entry $L_{ii}$ 
is the valency of vertex~$i$
and whose $(i,j)$-entry $L_{ij}$ is the negative of the number of
edges between vertices $i$ and $j$ if $i\ne j$.  
Then $L_{ii} = \sum_{j\ne i} L_{ij}$ for $1\leq i \leq n$, and so
the row sums of $\mathbf{L}$
are all zero.  It follows that $\mathbf{L}$ has eigenvalue~$0$
on the all-$1$ vector;  this is called the \emph{trivial eigenvalue} 
of~$\mathbf{L}$. We show below that the multiplicity of the zero eigenvalue
is equal to the number of connected components of $H$.  Thus the multiplicity
is~$1$ if and only if $H$ is connected.

Call the remaining eigenvalues of $\mathbf{L}$ \emph{non-trivial}.
They are all non-negative, as we show in the following theorem 
(see~\cite{bcc09}).

\begin{thm}
\label{thm:ledge}
\begin{enumerate}
\item
If~$\mathbf{L}$ is a Laplacian matrix, then $\mathbf{L}$~is 
positive semi-definite.
\item
If $\mathbf{L}$ is a Laplacian matrix of order $n$ and $\mathbf{x}$ is 
any vector in $\R^n$, then
\[
\mathbf{x}^\top \mathbf{L x} = \sum_{\mathrm{edges\ }ij} (x_i - x_j)^2.
\]
\item
If $\mathbf{L}_1$ and $\mathbf{L}_2$ are the Laplacian matrices of 
graphs $H_1$ and $H_2$ with the same vertices, and if $H_2$ is
obtained from $H_1$ by inserting one extra edge, then 
$\mathbf{L}_2- \mathbf{L}_1$ is positive semi-definite.
\item
If $\mathbf{L}$ is the Laplacian matrix of the graph $H$, then the 
multiplicity of the zero eigenvalue of $\mathbf{L}$ is equal to the 
number of connected components of~$H$.
\end{enumerate}
\end{thm}

\begin{pf}
Each edge between vertices $i$ and $j$ defines a $v\times v$ matrix whose 
entries are all $0$ apart from the following submatrix:
\[
\begin{array}{cc}
& \begin{array}{cc}
\hphantom{-}i & \hphantom{-}j
\end{array}\\
\begin{array}{c}i\\j\end{array}
&
\left[
\begin{array}{rr}
1 & -1\\
-1 & 1\rlap{\qquad .}
\end{array}
\right]
\end{array}
\]
The Laplacian is the sum of these matrices, which are all positive 
semi-definite.  This proves (a), (b) and~(c).

From (b), the vector $\mathbf{x}$ is in the null space of the
Laplacian if and only if $\mathbf{x}$ takes the same value on both 
vertices of each edge, which happens if 
and only if it takes a constant value on each connected component.  
This proves~(d).
\halmos
\end{pf}

Theorem~\ref{thm:ledge} shows that the smallest non-trivial
eigenvalue of a connected
graph is positive.  This eigenvalue is sometimes called the 
\textit{algebraic connectivity} of the graph.  The statistical
importance of this is shown in Section~\ref{sec:crit}.

In Section~\ref{sec:estvar} 
we shall need the Moore--Penrose generalized inverse of
$\mathbf{L}^-$  of $\mathbf{L}$
(see \cite{geninv}). 
Put $\mathbf{P}_0 = n^{-1}\mathbf{J}_n$, where $\mathbf{J}_n$ is the 
$n\times n$ matrix whose entries are all 
$1$, so that $\mathbf{P}_0$ is the matrix of orthogonal projection onto the
space spanned by the all-$1$ vector.  If $H$ is connected then 
$\mathbf{L}+\mathbf{P}_0$ is invertible, and 
\[
\mathbf{L}^- = \left(\mathbf{L}+\mathbf{P}_0 \right)^{-1} -\mathbf{P}_0,
\]
so that
$\mathbf{L}\mathbf{L}^- = \mathbf{L}^-\mathbf{L} = 
\mathbf{I}_n -\mathbf{P}_0$, where $\mathbf{I}_n$ is 
the identity matrix of order~$n$.

\subsection{Laplacians of the concurrence and Levi graphs}

There is a relationship between the Laplacian matrices of the concurrence
and Levi graphs of a block design $\Delta$. Let $\mathbf{N}$ be the
incidence matrix of the design, and $\mathbf{R}$ the diagonal matrix
(with rows and columns indexed by treatments) whose $(i,i)$ entry is the
replication $r_i$ of  treatment $i$. If the design is equireplicate,
then $\mathbf{R}=r\mathbf{I}_v$, where $r$ is the replication number.

For the remainder of the paper, we will use $\mathbf{L}$ for the Laplacian
matrix of the concurrence graph $G$ of $\Delta$, and $\tilde{\mathbf{L}}$ for
the Laplacian matrix of the Levi graph $\tilde{G}$ of $\Delta$. Then it is
straightforward to show that
\[\mathbf{L}=k\mathbf{R}-\mathbf{NN}^\top,\qquad
\tilde{\mathbf{L}}=\left[
\matrix{\mathbf{R} & -\mathbf{N}\cr
-\mathbf{N}^\top & k\mathbf{I}\cr}
\right]
.
\]

The Levi graph is connected if and only if the concurrence graph is 
connected; thus $0$ is a simple eigenvalue of $\tilde\mathbf{L}$ if and only
if it is a simple eigenvalue of $\mathbf{L}$,
which in turn occurs if and only if 
all contrasts between treatment parameters are estimable 
(see Section~\ref{sec:estvar}). A block design with this property is itself 
called \emph{connected}: we consider only connected block designs.

In the equireplicate case, the above expressions for $\mathbf{L}$ and
$\tilde\mathbf{L}$ give
a relationship between their
Laplacian eigenvalues, as follows. 
Let $\mathbf{x}$ be an eigenvector of $\mathbf{L}$
with eigenvalue $\phi\ne rk$. Then, for each of the two solutions $\theta$
of the quadratic equation
\[rk-\phi = (r-\theta)(k-\theta),\]
there is a unique vector $\mathbf{z}$ in $\R^b$
such that 
$[\begin{array}{cc}
\mathbf{x}^\top & \mathbf{z}^\top
\end{array}]^\top$
is an
eigenvector of $\tilde{\mathbf{L}}$ with eigenvalue $\theta$. Conversely,
any eigenvalue $\theta\ne k$ of $\tilde{\mathbf{L}}$ arises in this way.


The Laplacian matrices of the concurrence graphs in Fig.~\ref{fig:cgbin}
are shown in Table~\ref{tab:binary}.



\section{Statistical issues}
\subsection{Estimation and variance}
\label{sec:estvar}

As part of the experiment, we measure the response $Y_\omega$ on each 
experimental unit~$\omega$.  If $\omega$ is in block $\Gamma$, then we assume
that
\begin{equation}
Y_\omega = \tau_{f(\omega)} + \beta_\Gamma + \varepsilon_\omega;
\label{eq:linmod}
\end{equation}
here, $\tau_i$ is a constant depending on treatment~$i$, $\beta_\Gamma$ is a 
constant depending on block~$\Gamma$, and $\varepsilon_\omega$ is a random 
variable with expectation $0$ and variance~$\sigma^2$.  
Furthermore, if $\alpha\ne\omega$, 
then $\varepsilon_\alpha$ and $\varepsilon_\omega$ are uncorrelated.

It is clear that we can add a constant to every block parameter, and subtract
that constant from every treatment parameter, without 
changing~(\ref{eq:linmod}).  It is therefore impossible to estimate the
individual treatment parameters.  However, if the design is connected,
then we can estimate all \emph{contrasts} in the treatment parameters: 
that is, all linear combinations of the form $\sum_i x_i \tau_i$ for which
$\sum_i x_i = 0$.  In particular, we can estimate all the simple treatment
differences $\tau_i - \tau_j$.

An \emph{estimator} is a function of the responses $Y_\omega$, so it is itself
a random variable.  An estimator of a value is \emph{unbiased} if its 
expectation is equal to the true value; it is \emph{linear} if it is a linear
function of the responses.  Amongst linear unbiased estimators, the 
\emph{best} one (the so-called BLUE), is the one with the least variance.
Let $V_{ij}$ be the variance of the BLUE for $\tau_i - \tau_j$.

If all the experimental units form a single block, then the BLUE of 
${\tau_1 - \tau_2}$ is just the difference between the average responses for 
treatments $1$ and~$2$.
It follows that 
\[
V_{12} =\left(\frac{1}{r_1} + \frac{1}{r_2}\right)\sigma^2.
\] 
When $v=2$, this variance 
is minimized (for a given number of experimental units) when $r_1=r_2$.  
Moreover, if the responses are normally distributed then the length of the 
$95$\% confidence interval for $\tau_1 - \tau_2$ is proportional to 
$\mathrm{t}(r_1 + r_2 -2,0.975)\sqrt{V_{12}}$,
where $\mathrm{t}(d,p)$ is the $100p$-th percentile of the 
$\mathrm{t}$~distribution on $d$ degrees of freedom.  The smaller 
the confidence interval, 
the more likely is our estimate to be close to the true value.  This length
can be made smaller by 
increasing $r_1 + r_2$, 
decreasing $\card{r_1 - r_2}$, or 
decreasing $\sigma^2$.

However, matters are not so simple when $k<v$ and $v>2$.  The following result
can be found in any statistical textbook about block designs
(see the section on further reading for recommendations).

\begin{thm}
\label{thm:var}
Let $\mathbf{L}$ be the Laplacian matrix of the concurrence graph 
of a connected block design.  If $\sum_i x_i=0$, then the variance of 
the BLUE of $\sum_ix_i\tau_i$ is equal to 
$(\mathbf{x}^\top \mathbf{L}^{-} \mathbf{x})k\sigma^2$.
In particular, the variance~$V_{ij} $ of the BLUE of the simple
difference $\tau_i - \tau_j$ is given by
$V_{ij} = \left(L_{ii}^- +L_{jj}^{-} -2L_{ij}^{-}\right)k\sigma^2$.
\end{thm}

\subsection{Optimality criteria}
\label{sec:crit}

We want all of the $V_{ij}$ to be as small as possible, but this is a 
multi-dimensional problem if $v>2$.  Let $\bar V$ be the average of the 
variances $V_{ij}$ over all treatments $i$, $j$ with $i\ne j$.  
Theorem~\ref{thm:var} shows that, for each fixed~$i$,
\begin{eqnarray*}
\sum_{j\ne i} V_{ij} &=& 
\sum_{j\ne i} (L^-_{ii} + L^-_{jj} - 2L^-_{ij}) k\sigma^2\\
& = & [(v-1)L_{ii}^- + (\Tr(\mathbf{L}^-) - L^-_{ii}) + 2L_{ii}^-]k\sigma^2\\
& = & [vL^-_{ii} + \Tr(\mathbf{L}^{-})]k\sigma^2,
\end{eqnarray*}
because the row sums and column sums of $L$ are all $0$.  It follows that 
$\bar V = 
2k\sigma^2 \Tr(\mathbf{L}^-)/(v-1)$.

Let $\theta_1$, \ldots, $\theta_{v-1}$ be the non-trivial eigenvalues
of $\mathbf{L}$,  
now listed according to  multiplicity and in non-decreasing order.  Then
\[
\Tr(\mathbf{L}^-) = \frac{1}{\theta_1} + \cdots + \frac{1}{\theta_{v-1}},
\]
and so
\[
\bar V = 2k\sigma^2 \times 
\frac{1}{\mbox{harmonic mean of }\theta_1, \dots, \theta_{v-1}}.
\]

A block design is defined to be \emph{A-optimal} (in some given class of 
designs with the same values of $b$, $k$ and $v$) if it minimizes the value
of $\bar V$; here `A' stands for `average'.  Thus a design is A-optimal if 
and only if it maximizes the harmonic mean of $\theta_1$, \ldots, 
$\theta_{v-1}$. 

For $v>2$, the generalization of a confidence interval is a confidence 
ellipsoid centered at the point $(\hat \tau_1, \ldots, \hat \tau_v)$ 
which gives the estimated value of $(\tau_1, \ldots, \tau_v)$
in the $(v-1)$-dimensional subspace of $\R^v$ for which $\sum \tau_i = 0$.
A block design is called \emph{D-optimal} if it minimizes the volume of
this confidence ellipsoid.  Since this volume is proportional to 
$\sqrt{\det(\mathbf{L}^- + \mathbf{P}_0)}$, a design is D-optimal if 
and only if it maximizes the
geometric mean of $\theta_1$, \ldots, $\theta_{v-1}$. 
Here `D' stands for `determinant'.

Rather than looking at averages, we might consider the worst case.  If
all the entries in the vector $\mathbf{x}$ are multiplied by a 
constant~$c$, then
the variance of the estimator of $\sum x_i\tau_i$ is multiplied by $c^2$. 
Thus, those contrast vectors $\mathbf{x}$ which give the largest variance 
relative to their own length are those which maximize 
$\mathbf{x}^\top \mathbf{L}^- \mathbf{x}/\mathbf{x}^\top \mathbf{x}$;  
these are precisely the eigenvectors of $\mathbf{L}$ with 
eigenvalue $\theta_1$.
A design is defined to be \emph{E-optimal} if it maximizes the value of 
$\theta_1$; here `E' stands for `extreme'.

More generally, for $p$ in $(0, \infty)$, a design is called
\emph{$\Phi_p$-optimal} if it minimizes
\[
\left(\frac{\sum_{i=1}^{v-1} \theta_i^{-p}}{v-1}\right)^{1/p}.
\]
Thus A-optimality corresponds to $p=1$,
D-optimality corresponds to  the limit as $p\rightarrow 0$, and
E-optimality corresponds to the limit as $p\rightarrow \infty$.

Let $\mathbf{L}_1$ and $\mathbf{L}_2$ be the Laplacian matrices of 
the concurrence graphs of block designs $\Delta_1$ and $\Delta_2$ for 
$v$ treatments in blocks of size~$k$.  If $\mathbf{L}_2 -
\mathbf{L}_1$ is positive semi-definite, then $\Delta_2$ is at least
as good as $\Delta_1$ on all the $\Phi_p$-criteria.  
Theorem~\ref{thm:ledge}(c) shows that adding an extra block to a design 
cannot decrease its performance on any $\Phi_p$-criterion.

There are even more general classes of optimality criteria 
(see \cite{harman} and \cite{ss} for details).  
Here we concentrate on A-, D- and E-optimality.

\subsection{Questions and an example}
A first obvious question to ask is: do these criteria agree with each other?

Our optimality properties are all functions of the concurrence graph.
What features of this graph should we look for if we are searching for
optimal, or near-optimal, designs?
 Symmetry?
 (Nearly) equal degrees?
 (Nearly) equal numbers of edges between pairs of vertices?
 Distance-regularity?
 Large girth (ignoring cycles within a block)?
 Small numbers of short cycles (ditto)?
 High connectivity?
 Non-trivial automorphism group?

Is it more useful to look at the Levi graph rather than the concurrence graph?

\begin{eg}
\label{eg:cube}

Fig.~\ref{fig:cube} shows the values of the A- and D-criteria for all 
equi\-replicate block designs with $v=8$, $b=12$ and $k=2$: of course, these
are just regular graphs with $8$ vertices and degree~$3$. The harmonic mean
is shown on the $A$-axis, and the geometric mean on the $D$-axis.  
(Note that this figure includes some designs that were omitted from Figure~3 
of \cite{rabpaper}.)
The rankings on these two criteria are not exactly the same, but they do agree
at the top end, where it matters.  The second-best graph on both criteria is 
the cube; the best is the M\"obius ladder, whose vertices are the elements 
of $\Z_8$ and whose  edges are 
$\setof{i,i+1}$ and $\setof{i,i+4}$ for $i$ in $\Z_8$. These two graphs are so
close on both criteria that, for practical purposes, they can be regarded as
equally good.

\begin{figure}
\begin{center}
\setlength{\unitlength}{2cm}
\begin{picture}(5.4,3.2)(0.4,2.6)
\put(1,3){\line(1,0){4.5}}
\put(1,3){\line(0,1){2.5}}
\xlabel{1}{0.6}
\xlabel{2}{1.2}
\xlabel{3}{1.8}
\xlabel{4}{2.4}
\xlabel{5}{3.0}
\ylabel{3}{1.8}
\ylabel{4}{2.4}
\ylabel{5}{3.0}
\put(5.7,2.8){\makebox(0,0){$A$}}
\put(0.8,5.7){\makebox(0,0){$D$}}
\put(4.83,5.25){\makebox(0,0){$\boldsymbol{\times}$}}
\put(2.98,4.31){\makebox(0,0){$\boldsymbol{+}$}}
\put(3.35,4.66){\makebox(0,0){$\boldsymbol{+}$}}
\put(3.32,4.51){\makebox(0,0){$\boldsymbol{+}$}}
\put(3.89,4.95){\makebox(0,0){$\boldsymbol{+}$}}
\put(4.88,5.26){\makebox(0,0){$\boldsymbol{\times}$}}
\put(4.54,5.15){\makebox(0,0){$\boldsymbol{\times}$}}
\put(4.70,5.21){\makebox(0,0){$\boldsymbol{\times}$}}
\put(1.75,3.92){\makebox(0,0){$\boldsymbol{\circ}$}}
\put(3.60,4.81){\makebox(0,0){$\boldsymbol{+}$}}
\put(3.43,4.63){\makebox(0,0){$\boldsymbol{+}$}}
\put(4.01,4.91){\makebox(0,0){$\boldsymbol{+}$}}
\put(2.31,4.30){\makebox(0,0){$\boldsymbol{\circ}$}}
\put(4.08,4.95){\makebox(0,0){$\boldsymbol{+}$}}
\put(2.30,4.20){\makebox(0,0){$\boldsymbol{\circ}$}}
\put(2.44,4.45){\makebox(0,0){$\boldsymbol{\circ}$}}
\put(3.71,4.77){\makebox(0,0){$\boldsymbol{+}$}}
\put(4.19,4.99){\makebox(0,0){$\boldsymbol{+}$}}
\put(4.9,5.4){\makebox(0,0)[l]{M\"obius ladder}}
\put(4.8,5.1){\makebox(0,0)[l]{cube}}
\put(4.18,4.75){\circle*{0.07}}
\put(4.3,4.7){\makebox(0,0)[l]{$K_{2,6}$}}
\end{picture}
\end{center}
\caption{Values of two optimality criteria for all equireplicate block 
designs with $v=8$, $b=12$, and $k=2$, and for $K_{2,6}$}
\label{fig:cube}
\end{figure}
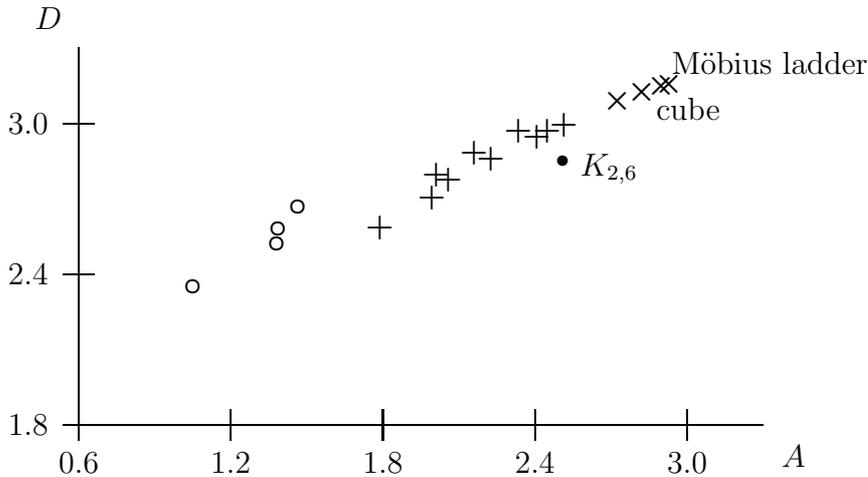

The plotting symbols show the edge-connectivity of the graphs: 
edge-connectivity $3$, $2$, $1$ is shown as $\boldsymbol{\times}$, 
$\boldsymbol{+}$, $\boldsymbol{\circ}$ respectively.  This does suggest
that the higher the edge-connectivity the better is the design on the
A- and D- criteria.  This is intuitively reasonable: if $k=2$, then the 
edge-connectivity is the minimum number of blocks whose removal disconnects
the design.  In this context, it has been called \emph{breakdown number}: 
see \cite{mahbub}.

The four graphs with edge-connectivity $3$ have no double edges, so 
concurrences differ by at most~$1$.  The only other regular graph with no
double edges is ranked eighth (amongst regular graphs) by the A-criterion.
This suggests that (near-)equality of concurrences is not sufficient to
give a good design.

The symbol 
{\setlength{\unitlength}{2cm}\raisebox{1mm}{\circle*{0.07}}}
shows the non-regular graph $K_{2,6}$, which also has eight vertices and twelve
edges.  It is not as good as the regular graphs with edge-connectivity~$3$, 
but it beats many of the other regular graphs.

This pattern is typical of the block designs investigated by statisticians
for most of the 20th century. The A- and D-criteria agree closely at the
top end.  High edge-connectivity appears to show good designs.  Many of
the best designs have a high degree of symmetry.
\end{eg}

\section{Highly patterned block designs}

\subsection{Balanced incomplete-block designs}
\label{sec:bibd}

BIBDs are intuitively appealing, as they seem to give equal weight to all
treatment comparisons.  They were introduced for agricultural experiments
by Yates in \cite{fy:ibd}.

In \cite{ksh}, Kshirsagar proved that, if there exists a BIBD for given
values of $v$, $b$ and $k$, then it is A-optimal.  Kiefer generalized this
in \cite{kiefer} to cover $\Phi_p$-optimality for all $p$ in $(0,\infty)$, 
including the limiting cases of D- and E-optimality.  The core of Kiefer's 
proof is as follows:  binary designs maximize $\Tr(\mathbf{L})$, 
which is equal to $\sum_{i=1}^{v-1} \theta_i$;  
for any fixed value $T$ of this sum of positive
numbers, $\sum \theta_i^{-p}$ is minimized at $[T/(v-1)]^{-p}$ 
when $\theta_1 = \cdots = \theta_{v-1} = T/(v-1)$; and $T^{-p}$ is minimized 
when $T$ is maximized.

\subsection{Other special designs}

Of course, it frequently occurs that the values of $b,v,k$ available for
an experiment are such that no BIBD exists. (Necessary conditions for the
exist\-ence of a BIBD include the well-known divisibility conditions 
$v\mid bk$ and $v(v-1)\mid bk(k-1)$, which follow from the elementary 
results in Section~\ref{sec:intro}, 
and \emph{Fisher's inequality} asserting that $b\ge v$.)

In the absence of a BIBD, various other special types of design have been
considered, and some of these have been proved optimal. Here is a short
sample.

A design is \emph{group-divisible} if the treatments can be partitioned into
``groups'' all of the same size, so that the number of blocks containing two
treatments is $\lambda_1$ if they belong to the same group and $\lambda_2$
otherwise. Ch\^eng~\cite{cheng1,cheng2} showed that if there is a
group-divisible design with two groups and $\lambda_2=\lambda_1+1$ 
in the class of designs with given $v$, $b$, $k$, then it is 
$\Phi_p$-optimal for all $p$, and in particular is
A-, D- and E-optimal.

A \emph{regular-graph design} is a binary equireplicate design with two
possible concurrences $\lambda$ and $\lambda+1$. It is easily proved that,
in such a design, the number of treatments lying in $\lambda+1$ blocks with
a given treatment is constant; so the graph $H$ whose vertices are the 
treatments, two vertices joined if they lie in $\lambda+1$ blocks, 
is regular.

Now Ch\^eng~\cite{cheng2} showed that a group-divisible design with
$\lambda_2 = \lambda_1 + 1 = 1$, if one
exists, is $\Phi_p$-optimal in the class of regular-graph designs
for all~$p$.  Cheng and
Bailey~\cite{cheng_bailey} showed that a regular-graph design for which the
graph is \emph{strongly regular} (see~\cite{pjc:sr}) and which has singular
concurrence matrix is $\Phi_p$-optimal, for all~$p$, among binary
equireplicate designs with given $v$, $b$, $k$.

Designs with the property described here are particular examples of 
\emph{partially balanced designs} with respect to an association scheme:
see Bailey~\cite{rab:as}.

Another class which has turns out to be optimal in many cases, but whose
definition is less combinatorial, consists of the \emph{variance-balanced
designs}, which we consider later in the chapter.

\section{Graph concepts linked to D-optimality}
\label{sec:span}

\subsection{Spanning trees of the concurrence graph}

Let $G$ be the concurrence graph of a connected block design, and let 
$\mathbf{L}$ be its Laplacian matrix.  A \emph{spanning tree} for $G$ is a 
spanning subgraph which is a tree.
Kirchhoff's famous Matrix-tree theorem in \cite{old} 
states the following:

\begin{thm}
\label{thm:span}
If $G$ is a connected graph with $v$~vertices and Laplacian 
matrix~$\mathbf{L}$, then the product of the non-trivial eigenvalues 
of~$\mathbf{L}$ is equal to $v$ multiplied by the number of spanning 
trees for $G$.
\end{thm}

Thus we have a test for D-optimality:
\begin{quote}
\textit{A design is D-optimal if and only if its concurrence graph has the 
maximal number of spanning trees.}
\end{quote}

Note that Theorem~\ref{thm:span} gives an easy proof of Cayley's theorem 
on the number of spanning trees for the complete graph $K_v$.  The 
non-trivial eigenvalues of its Laplacian matrix are all equal to $v$, so
Theorem~\ref{thm:span} shows that it has $v^{v-2}$ spanning trees.

If $G$ is sparse, it may be much easier to count the number of spanning 
trees than to compute the eigenvalues of~$\mathbf{L}$.  For example, if 
$G$ has a single cycle, which has length~$s$, then the number of 
spanning trees is~$s$, irrespective of the remaining edges in $G$.

In the context of optimal block designs, Gaffke discovered the importance
of Kirchhoff's theorem in \cite{gaff}. Cheng followed this up in papers 
such as \cite{cheng1,cheng2,cheng3}.  Particularly intriguing is the 
following theorem from \cite{cmw}.

\begin{thm}
Consider block designs with $k=2$ (connected graphs).  For each 
given~$v$ there is a threshold $b_0$ such that if $b\geq b_0$ then 
any D-optimal design for $v$ treatments in $b$ blocks of size~$2$ 
is \emph{nearly balanced} in the sense that
\begin{itemize}
\item
no pair of replications differ by more than $1$;
\item
for each fixed $i$, no pair of concurrences $\lambda_{ij}$ differ
by more than $1$.
\end{itemize}
\end{thm}

In fact, there is no known example with $b_0 >v-1$, which is the minimal
number of blocks required for connectivity.

\subsection{Spanning trees of the Levi graph}

In \cite{gafLevi} Gaffke stated the following relationship between the
numbers of spanning trees in the concurrence graph and the Levi graph.

\begin{thm}
\label{thm:gaff}
Let $G$ and $\tilde G$ be the concurrence graph and Levi graph for a
connected incomplete-block design for $v$ treatments in $b$ blocks of
size $k$.  Then the number of spanning trees for $\tilde G$ is equal to
$k^{b-v+1}$ times the number of spanning trees for $G$.
\end{thm} 

Thus, an alternative test for D-optimality is to count the number of
spanning trees in the Levi graph.  For binary designs, the Levi graph
has fewer edges than the concurrence graph if and only if $k\geq 4$.

\section{Graph concepts linked to A-optimality}
\subsection{The concurrence graph as an electrical network}
\label{sec:concelec}

We can consider the concurrence graph $G$ as an electrical network with a
1-ohm resistance in each edge.  Connect a $1$-volt battery between
vertices $i$ and $j$.
Then current flows in the network, according to these rules.
\begin{description}
  \item[Ohm's Law:] 
In every edge, the voltage drop is the product of the current and the 
resistance.
\item[Kirchhoff's Voltage Law:]
 The total voltage drop from one vertex to any other vertex is the same 
no matter which path we take from one to the other.
\item[Kirchhoff's Current Law:] At each vertex which is not
 connected to the battery, the total current coming in is equal
to the total current going out.
\end{description}

We find the total current  from $i$ to $j$, and
then use Ohm's Law to define
the effective resistance $R_{ij}$ between $i$ and $j$ as the reciprocal of 
this current. It is a standard result of electrical network theory that the
linear equations implicitly defined above for the currents and voltage
differences have
a unique solution.

Let $\mathcal{T}$ be the set of treatments and $\Omega$ the set of
experimental units.  Current flows in each edge $e_{\alpha\omega}$,
where $\alpha$ and $\omega$ are experimental units in the same block which
receive different treatments; let $I(\alpha,\omega)$ be the current
from $f(\alpha)$ to $f(\omega)$ in this edge. Thus $I$ is a function 
$I \colon \Omega \times \Omega \mapsto \R$ such that 
\begin{enumerate}
  \item $I(\alpha,\omega) = 0$ if $\alpha=\omega$ or if $f(\alpha) =
    f(\omega)$   or if $\alpha$ and $\omega$ are in different blocks.
\item
$I(\alpha,\omega) = -I(\omega,\alpha)$ 
for $(\alpha,\omega)$ in $\Omega\times \Omega$. 
\end{enumerate}
This defines a further function $I_\mathrm{out} \colon \mathcal{T} \mapsto \R$ 
by
\[
I_\mathrm{out} (l) = \sum_{\alpha: f(\alpha)=l} \ \sum_{\omega\in\Omega} 
I(\alpha,\omega) \quad \mbox{for $l$ in $\mathcal{T}$}. 
\]
Voltage is another function $V\colon \mathcal{T} \mapsto \R$.
The following two conditions ensure that Ohm's and Kirchhoff's Laws are 
satisfied.
\begin{enumerate}
  \addtocounter{enumi}{2}
\item
If there is any edge in $G$ between $f(\alpha)$ and $f(\omega)$, then
\[I(\alpha,\omega) = V(f(\alpha)) - V(f(\omega)).\]
\item
If $l\notin\setof{i,j}$, then $I_{\mathrm{out}}(l)=0$.
\end{enumerate}

If $G$ is connected and different voltages $V(i)$ and $V(j)$ are given for 
a pair of distinct treatments $i$ and $j$, then there are unique functions 
$I$ and $V$ satisfying conditions (a)--(d). Moreover, $I_\mathrm{out}(j) 
= -I_{\mathrm{out}} (i) \ne 0$.  Then $R_{ij}$  is defined by
\[
R_{ij} = \frac{V(i) - V(j)}{I_{\mathrm{out}}(i)}.
\]
It can be shown that the value of $R_{ij}$ does not depend on the choice of
values for $V(i)$ and $V(j)$, so long as these are different.  In practical 
examples, it is usually convenient to take $V(i)=0$ and let $I$ take integer 
values.

What has all of this got to do with block designs?  The following theorem, 
which is a standard result from electrical engineering, gives the answer.

\begin{thm}
\label{thm:resist}
If $\mathbf{L}$ is the Laplacian matrix of a connected graph $G$, 
then the effective 
resistance $R_{ij}$ between vertices $i$ and $j$ is given by
\[R_{ij} = \left(L_{ii}^- +L_{jj}^{-} -2L_{ij}^{-}\right).\]
\end{thm}

Comparing this with Theorem~\ref{thm:var}, we see that
$V_{ij} = R_{ij} \times k\sigma^2$. Hence we have a test for A-optimality:

\begin{quote}
\textit{A design is A-optimal if and only if its concurrence graph, regarded
as an electrical network, minimizes the sum of the pairwise effective 
resistances between all pairs of vertices.}
\end{quote}

Effective resistances are easy to calculate without matrix inversion
if the graph is sparse.

\begin{figure}
\begin{center}
  \setlength{\unitlength}{5.8mm}
\begin{picture}(15,14)(-6,-7)
\put(0,6.67){\blobb}
\put(0,-6.67){\blobb}
\put(-6,-3.33){\blobb}
\put(-2,-3.33){\blobb}
\put(6,-3.33){\blobb}
\put(2,-3.33){\blobb}
\put(-6,3.33){\blobb}
\put(-2,3.33){\blobb}
\put(6,3.33){\blobb}
\put(2,3.33){\blobb}
\put(-4,0){\blobb}
\put(4,0){\blobb}
\put(-6,-3.33){\line(1,0){12}}
\put(-6,3.33){\line(1,0){12}}
\put(-6,-3.33){\line(3,5){6}}
\put(6,-3.33){\line(-3,5){6}}
\put(-6,3.33){\line(3,-5){6}}
\put(6,3.33){\line(-3,-5){6}}
\put(0,-3.33){\vector(1,0){0}}
\put(0,-2.9){\makebox(0,0){14}}
\put(4,-3.33){\vector(1,0){0}}
\put(4,-2.9){\makebox(0,0){7}}
\put(-4,-3.33){\vector(-1,0){0}}
\put(-4,-2.9){\makebox(0,0){5}}
\put(-1,-5){\vector(4,-3){0}}
\put(-1.4,-5.2){\makebox(0,0){7}}
\put(-3,-1.67){\vector(-4,3){0}}
\put(-2.6,-1.4){\makebox(0,0){10}}
\put(-5,1.67){\vector(-4,3){0}}
\put(-5.4,1.7){\makebox(0,0){5}}
\put(1,-5){\vector(4,3){0}}
\put(1.4,-5.2){\makebox(0,0){7}}
\put(3,-1.67){\vector(4,3){0}}
\put(2.6,-1.4){\makebox(0,0){14}}
\put(5,1.67){\vector(4,3){0}}
\put(5.4,1.7){\makebox(0,0){19}}
\put(-5,-1.67){\vector(4,3){0}}
\put(-5.4,-1.4){\makebox(0,0){5}}
\put(-3,1.67){\vector(4,3){0}}
\put(-2.6,1.4){\makebox(0,0){10}}
\put(-1,5){\vector(4,3){0}}
\put(-1.4,5.2){\makebox(0,0){5}}
\put(0,3.33){\vector(1,0){0}}
\put(0,2.9){\makebox(0,0){10}}
\put(4,3.33){\vector(1,0){0}}
\put(4,2.9){\makebox(0,0){17}}
\put(-4,3.33){\vector(1,0){0}}
\put(-4,2.9){\makebox(0,0){5}}
\put(5,-1.67){\vector(-4,3){0}}
\put(5.4,-1.4){\makebox(0,0){7}}
\put(3,1.67){\vector(-4,3){0}}
\put(2.6,1.4){\makebox(0,0){2}}
\put(1,5){\vector(4,-3){0}}
\put(1.4,5.2){\makebox(0,0){5}}
\put(-1.7,-2.83){\makebox(0,0){$i$}}
\put(-2.6,-4){\makebox(0,0){$[0]$}}
\put(2.6,-4){\makebox(0,0){$[-14]$}}
\put(-2.6,4){\makebox(0,0){$[-20]$}}
\put(2.6,4){\makebox(0,0){$[-30]$}}
\put(-5.1,0){\makebox(0,0){$[-10]$}}
\put(5.1,0){\makebox(0,0){$[-28]$}}
\put(5.4,2.9){\makebox(0,0){$j$}}
\put(6.6,4){\makebox(0,0){$[-47]$}}
\put(0,-7.33){\makebox(0,0){$[-7]$}}
\put(0,7.33){\makebox(0,0){$[-25]$}}
\put(-6.6,-4){\makebox(0,0){$[-5]$}}
\put(-6.6,4){\makebox(0,0){$[-15]$}}
\put(6.6,-4){\makebox(0,0){$[-21]$}}
\end{picture}
\end{center}
\caption{
The current between vertices $i$ and $j$ in a 
concurrence graph}
\label{fig:star}
\end{figure}

Figure~\ref{fig:star} shows the concurrence graph of a block design with
$v=12$, $b=6$ and $k=3$.  Only vertices $i$ and $j$ are labelled.  Otherwise,
numbers beside arrows denote current and numbers in square brackets denote
voltage.  It is straightforward to check that conditions (a)--(d) are 
satisfied.  Now $V(i) - V(j) = 47$ and $I_{\mathrm{out}}(i) = 36$, 
and so $R_{ij} = 47/36$.  Therefore $V_{ij} = (47/12)\sigma^2$.
Moreover, for graphs consisting of $b$~triangles arranged in a cycle 
like this, it is clear that average effective resistance, and hence
the average pairwise variance, can be calculated as a function of $b$. 

\subsection{The Levi graph as an electrical network}
\label{sec:levielec}

The Levi graph $\tilde G$ of a block design can also be considered as an 
electrical network.  Denote by $\mathcal{B}$ the set of blocks.  Now current
is defined on the ordered edges of the Levi graph. Recall that, if $\omega$
is an experimental unit in block $\Gamma$, then the edge $\tilde{e}_{\omega}$
joins $\Gamma$ to $f(\omega)$. Thus current is defined on 
$(\Omega \times \mathcal{B}) \cup (\mathcal{B} \times \Omega)$ 
and voltage is defined on $\mathcal{T} \cup \mathcal{B}$.  
Conditions (a)--(d) in Section~\ref{sec:concelec} need to be modified 
appropriately.

The next theorem shows that a current--voltage pair $(I,V)$ on the concurrence
graph $G$ can be transformed into a current--voltage pair 
$(\tilde I, \tilde V)$ on the Levi graph $\tilde G$.  In $\tilde G$, the
current $\tilde I(\alpha,\Gamma)$ flows in edge $\tilde{e}_\alpha$ 
from vertex $f(\alpha)$ to vertex $\Gamma$, where $\alpha\in\Gamma$. 
Hence the pairwise variance $V_{ij}$ can also be
calculated from the effective resistance $\tilde R_{ij}$ in the Levi graph.

\begin{thm}
\label{thm:tjur}
Let $G$ 
be the concurrence graph and 
$\tilde G$ be the Levi graph of a connected 
block design with block size~$k$.  If\/ $i$ and $j$ are two distinct 
treatments, let $R_{ij}$ and $\tilde R_{ij}$ be the effective resistance 
between vertices $i$ and $j$ in the electrical networks defined by $G$ and 
$\tilde G$, respectively.  Then $\tilde R_{ij} = kR_{ij}$, and so 
$V_{ij} = \tilde R_{ij} \sigma^2$.
\end{thm}

\begin{pf}
Let $(I,V)$ be a current--voltage pair on $G$.  
For $(\alpha,\Gamma) \in \Omega \times \mathcal{B}$, put
\[
\tilde I(\alpha,\Gamma) = -\tilde I(\Gamma,\alpha) = 
\sum_{\omega\in\Gamma} I(\alpha,\omega)
\]
if $\alpha\in\Gamma$; otherwise, put
$\tilde I(\alpha,\Gamma) = \tilde I(\Gamma,\alpha) = 0$.
Put $\tilde V(i) = kV(i)$ for all $i$ in~$\mathcal{T}$, and
\[
\tilde V(\Gamma) = \sum_{\omega\in\Gamma} V(f(\omega))
\]
for all $\Gamma$ in $\mathcal{B}$. 
It is clear that $\tilde I$ satisfies the analogues of conditions (a) and~(b).

If $\alpha\in\Gamma$, then
\begin{eqnarray*}
\tilde I(\alpha,\Gamma) &=& 
\sum_{\omega\in\Gamma} I(\alpha,\omega)\\
& = & \sum_{\omega\in\Gamma}[V(f(\alpha)) - V(f(\omega))]\\
& = & kV(f(\alpha)) - \tilde V(\Gamma) 
= \tilde V(f(\alpha)) - \tilde V(\Gamma),
\end{eqnarray*}
so the analogue of condition (c) is satisfied.

If $\Gamma\in\mathcal{B}$, then
\[
\tilde I_{\mathrm{out}}(\Gamma) = 
\sum_{\alpha\in\Gamma} \tilde I(\Gamma,\alpha)
= -\sum_{\alpha\in\Gamma} \sum_{\omega\in\Gamma} I(\alpha,\omega) = 0,
\]
because $I(\alpha,\alpha)=0$ and $I(\alpha,\omega) = -I(\omega,\alpha)$.
If $l\in \mathcal{T}$ then
\[
\tilde I_{\mathrm{out}}(l) =  
\sum_{\alpha: f(\alpha)=l} \ \sum_{\Gamma\in\mathcal{B}} 
\tilde I(\alpha,\Gamma) =
\sum_{\alpha: f(\alpha)=l}\ \sum_{\omega\in\Omega} I(\alpha,\omega) = 
I_{\mathrm{out}}(l).
\]
In particular, $\tilde I_{\mathrm{out}}(l) =  0$ if 
$l\notin\setof{i,j}$, which shows that the analogue of condition (d)
is satisfied. It follows that $(\tilde I,\tilde V)$ is the current--voltage
pair on $\tilde G$ defined by $\tilde V(i)$ and $\tilde V(j)$.

Now
\[
\tilde R_{ij} = \frac{\tilde V(i) - \tilde V(j)}{\tilde I_{\mathrm{out}}(i)}
= \frac{k(V(i)-V(j))}{I_{\mathrm{out}}(i)} = kR_{ij}.
\]
Then Theorems~\ref{thm:var} and~\ref{thm:resist} show that 
$V_{ij} = R_{ij}\sigma^2$.
\halmos
\end{pf}

When $k=2$ it seems to be easier to use the concurrence graph than
the Levi graph, because it has
fewer vertices, but for larger values of $k$ the Levi graph may be better, as 
it does not have all the within-block cycles that the concurrence graph has.
Fig.~\ref{fig:levistar} gives the Levi graph of the block design whose 
concurrence graph is in Fig.~\ref{fig:star}, with the same two vertices 
$i$ and $j$ attached to the battery.  This gives $\tilde R_{ij} = 47/12$, which
is in accordance with Theorem~\ref{thm:tjur}.

\begin{figure}[htbp]
\begin{center}
\setlength{\unitlength}{8mm}
\begin{picture}(10,6)(-0.5,0)
\multiput(0,1)(2,0){6}{\blobb}
\multiput(0,3)(2,0){6}{\blobb}
\multiput(0,5)(2,0){6}{\blobb}
\put(-0.5,1){\makebox(0,0)[r]{treatments}}
\put(-0.5,3){\makebox(0,0)[r]{blocks}}
\put(-0.5,5){\makebox(0,0)[r]{treatments}}
\multiput(0,1)(2,0){6}{\line(0,1){4}}
\put(0,1){\line(5,1){10}}
\multiput(2,1)(2,0){5}{\line(-1,1){2}}
\multiput(2,2.1)(2,0){3}{\vector(0,1){0}}
\put(6,4.4){\vector(0,1){0}}
\put(0,1.9){\vector(0,-1){0}}
\put(8,1.9){\vector(0,-1){0}}
\put(10,1.9){\vector(0,-1){0}}
\put(3.1,1.9){\vector(1,-1){0}}
\put(4.7,2.3){\vector(1,-1){0}}
\put(0.9,2.1){\vector(-1,1){0}}
\put(6.9,2.1){\vector(-1,1){0}}
\put(8.9,2.1){\vector(-1,1){0}}
\put(4.7,1.95){\vector(4,1){0}}
\put(0,0.7){\makebox(0,0)[t]{$[-10]$}}
\put(2,0.7){\makebox(0,0)[t]{$[0]$}}
\put(2,0){\makebox(0,0)[t]{$i$}}
\put(4,0.7){\makebox(0,0)[t]{$[-14]$}}
\put(6,0.7){\makebox(0,0)[t]{$[-28]$}}
\put(8,0.7){\makebox(0,0)[t]{$[-30]$}}
\put(10,0.7){\makebox(0,0)[t]{$[-20]$}}
\put(0.7,3.3){\makebox(0,0){$[-5]$}}
\put(2.7,3.3){\makebox(0,0){$[-7]$}}
\put(4.7,3.3){\makebox(0,0){$[-21]$}}
\put(6.7,3.3){\makebox(0,0){$[-35]$}}
\put(8.7,3.3){\makebox(0,0){$[-25]$}}
\put(10.7,3.3){\makebox(0,0){$[-15]$}}
\put(0,5.3){\makebox(0,0)[b]{$[-5]$}}
\put(2,5.3){\makebox(0,0)[b]{$[-7]$}}
\put(4,5.3){\makebox(0,0)[b]{$[-21]$}}
\put(6,5.3){\makebox(0,0)[b]{$[-47]$}}
\put(6,6){\makebox(0,0)[b]{$j$}}
\put(8,5.3){\makebox(0,0)[b]{$[-25]$}}
\put(10,5.3){\makebox(0,0)[b]{$[-15]$}}
\put(-0.3,2){\makebox(0,0){5}}
\put(0.7,2){\makebox(0,0){5}}
\put(1.7,2){\makebox(0,0){7}}
\put(2.7,2){\makebox(0,0){7}}
\put(3.7,2){\makebox(0,0){7}}
\put(4.6,1.6){\makebox(0,0){5}}
\put(4.8,2.5){\makebox(0,0){7}}
\put(6.3,2){\makebox(0,0){7}}
\put(7.3,2){\makebox(0,0){5}}
\put(8.3,2){\makebox(0,0){5}}
\put(9.3,2){\makebox(0,0){5}}
\put(10.3,2){\makebox(0,0){5}}
\put(5.5,4.3){\makebox(0,0){12}}
\end{picture}
\end{center}
\caption{Current between $i$ and $j$ for the Levi graph corresponding to
the concurrence graph in Fig.~\ref{fig:star}}
\label{fig:levistar}
\end{figure}

Here is another way of visualizing Theorem~\ref{thm:tjur}.  From the block
design we construct a graph $G_0$ with vertex-set $\mathcal{T} \cup \Omega 
\cup\mathcal{B}$: the edges are $\setof{\alpha,\Gamma}$ for $\alpha\in\Gamma 
\in \mathcal{B}$ and $\setof{\alpha,f(\alpha)}$ for $\alpha\in\Omega$.
Let $(I_0,V_0)$ be a current--voltage pair on $G_0$ for which both battery
vertices are in $\mathcal{T}$.  We obtain the Levi graph $\tilde G$ from $G_0$ 
by becoming blind to the vertices in~$\Omega$.  Thus the resistance in each 
edge of $\tilde G$ is twice that in each edge in $G_0$, so this step
multiplies each effective resistance by $2$.  

Because none of the battery
vertices is in~$\mathcal{B}$, we can now obtain $G$ from $\tilde G$ by 
replacing each path of the form $(i,\Gamma,j)$ by an edge $(i,j)$.
There is no harm in scaling all the voltages by the same amount, so we
can obtain $(I,V)$ on $G$ from $(\tilde I, \tilde V)$ on $\tilde G$ by
putting $V(i) = \tilde V(i)/k$ for $i$ in $\mathcal{T}$, and
$I(\alpha,\omega) = V(f(\alpha)) - V(f(\omega))$ for $\alpha$, $\omega$
in the same block.  If $\Gamma$ is a block, then
\[
0 = \sum_{\alpha\in\Gamma} \tilde I(\alpha,\Gamma) 
= \sum_{\alpha\in\Gamma}[\tilde V(f(\alpha)) - \tilde V(\Gamma)]
= k\sum_{\alpha\in\Gamma} V(f(\alpha)) - k\tilde V(\Gamma),
\]
and so $\tilde V(\Gamma) = \sum_{\alpha\in\Gamma} V(f(\alpha))$.
Also, if $\alpha\in\Gamma$, then
\begin{eqnarray*}
\sum_{\omega\in\Gamma} I(\alpha,\omega)& =&
\sum_{\omega\in\Gamma}[V(f(\alpha)) - V(f(\omega))]\\
&=& kV(f(\alpha)) - \sum_{\omega\in\Gamma}V(f(\omega)) \\
&=& \tilde V(f(\alpha)) - \tilde V (\Gamma)\\& =& \tilde I(\alpha,\Gamma).
\end{eqnarray*}
Therefore, this transformation reverses the one used in the proof of
Theorem~\ref{thm:tjur}. 

There is yet another way of obtaining Theorem~\ref{thm:tjur}.  If we use 
the responses $Y_\omega$ to estimate the block parameters $\beta_\Gamma$ in
(\ref{eq:linmod}) as well as the treatment parameters $\tau_i$, then standard
theory of linear models shows that, if the design is connected, then we
can estimate linear combinations of the form
$\sum_{i=1}^v x_i \tau_i + \sum_{j=1}^b z_j\beta_j$
so long as $\sum x_i = \sum z_j$.  Moreover, the variance of the BLUE of
this linear combination is
\[
[\begin{array}{cc}\mathbf{x}^\top & \mathbf{z}^\top\end{array}] \mathbf{C}^-
\left[ 
\begin{array}{c}\mathbf{x} \\ \mathbf{z}\end{array}
\right]
\sigma^2,
\qquad\mbox{where} \qquad
\mathbf{C} = 
\left[
\begin{array}{lc}
\mathbf{R} & \mathbf{N}\\
\mathbf{N}^\top & k\mathbf{I}_b
\end{array}
\right]
\]
and $\mathbf{R}$ is the diagonal matrix of replications.

If we reparametrize equation~(\ref{eq:linmod}) by replacing $\beta_j$ by
$-\gamma_j$ for $j=1$, \ldots, $b$, then the estimable quantities are the
contrasts in $\tau_1$, \ldots, $\tau_v$, $\gamma_1$, \ldots, $\gamma_b$.  The
so-called \emph{information matrix} $\mathbf{C}$ must be modified by
multiplying the  
last $b$~rows and the last $b$~columns by $-1$: this gives 
precisely the Laplacian $\tilde \mathbf{L}$ of the Levi
graph $\tilde G$. Just as for~$\mathbf{L}$, but unlike $\mathbf{C}$,
the null space is spanned by the all-$1$ vector. 

\subsection{Spanning thickets}

We have seen that the value of the D-criterion is a function of 
the number of spanning trees of the concurrence graph~$G$. It turns out
that the closely related  
notion of a spanning thicket enables us to calculate the A-criterion; more 
precisely, the value of each pairwise effective resistance in $G$.

A \emph{spanning thicket} for the graph is a spanning subgraph that consists
of two trees (one of them may be an isolated vertex).

\begin{thm}
\label{thm:thick}
If $i$ and $j$ are distinct vertices of $G$ then
\[
R_{ij} =
 \frac{\mbox{\normalfont{number of spanning thickets with $i$, $j$ in different
  parts}}}{\mbox{\normalfont{number of spanning trees}}}
\ .
\]
\end{thm}

This is also rather easy to calculate directly when the graph is sparse.

Summing all the $R_{ij}$ and using Theorem~\ref{thm:thick} gives the following 
result from \cite{shap}.

\begin{thm}
If $F$ is a spanning thicket for the concurrence graph~$G$, denote by
$F_1$ and $F_2$ the sets of vertices in its two trees. Then
\[
\sum_{i< j} R_{ij} =\frac{\displaystyle
\sum_{\mathrm{spanning\ thickets}\ F} \left| F_1 \right|
    \left| F_2\right|}{\mbox{\normalfont{number of spanning trees}}}
\]
\end{thm}

\subsection{Random walks and electrical networks}

It was first pointed out by Kakutani in 1945 that there is a very close
connection between random walks and electrical networks. In a simple random
walk, a single step works as follows: starting at a vertex, we choose an edge
containing the vertex at random, and move along it to the other end. This
definition accommodates multiple edges, and is easily adapted to
graphs with edge weights (where the probability of moving along an edge is
proportional to the weight of the edge).

If we are thinking of an edge-weighted graph as an electrical network, we
take the weights to be the conductances of the edges (the reciprocals of the
resistances).

The connection is simple to state:

\begin{thm}
Let $i$ and $j$ be distinct vertices of the connected edge-weighted graph $G$.
Apply voltages of $1$ at $i$ and $0$ at $j$. Then the voltage at a
vertex $l$ is equal to the probability that the random walk, starting at $l$,
reaches $i$ before it reaches $j$.
\end{thm}

From this theorem, it is possible to derive a formula for the effective
resistance between two vertices. Here are two such formulas. Given two
vertices $i$ and $j$, let $P_\mathrm{esc}(i\to j)$ be the probability that
a random walk starting at $i$ reaches $j$ before returning to $i$; and let
$S_i(i,j)$ be the expected number of times that a random walk starting at $i$
visits $i$ before reaching $j$. Then the effective resistance between $i$ and
$j$ is given by either of the two expressions 
\[
\frac{1}{d_iP_\mathrm{esc}(i\to j)}\qquad\mbox{and}\qquad 
\frac{S_i(i,j)}{d_i},
\]
where $d_i$ is the degree of $i$. (If the edge resistances
are not all $1$, then the term $d_i$ should be replaced by the sum of the
reciprocals of the resistances of all edges incident with vertex $i$.)

The random walk approach gives alternative
proofs of some of the main results about electrical networks. We discuss this
further in the guide to the literature.

\subsection{Foster's formula and generalizations}

In 1948, Foster~\cite{foster} discovered that the sum of the effective
resistances between all \emph{adjacent} pairs of vertices of a connected graph
on $v$ vertices is equal to $v-1$. Thirteen years later, he found a similar
formula for pairs of vertices at distance~$2$:
\[\sum_{i\sim h\sim j}\frac{R_{ij}}{d_h} = v-2.\]
Further extensions have been found, but require a stronger condition on the
graph. The sum of resistances between all pairs of vertices at distance at
most~$m$ can be written down explicitly if the graph is \emph{walk-regular
up to distance $m$}; this means that the number of closed walks of length $k$
starting and finishing at a vertex $i$ is independent of $i$, for $k\le m$.
The formula was discovered by Emil Vaughan, to whom this part of the chapter
owes a debt.

In particular, if the graph is distance-regular (see \cite{bcn}), 
then the value of the A-criterion can be written down in terms of the 
so-called \emph{intersection array} of the graph.

\subsection{Distance}

At first sight it seems obvious that pairwise variance should decrease
as concurrence increases, but there are many counter-examples to
this. However, the following theorem is proved in \cite{rab:as}.

\begin{thm}
If the Laplacian matrix $\mathbf{L}$ has precisely two distinct
non-trivial eigenvalues, then pairwise variance is a decreasing linear
function of concurrence.
\end{thm}

It does appear that effective resistance, and hence pairwise variance, 
generally increases with distance in the concurrence graph. 
In \cite[Question 5.1]{bcc09} we pointed out that this is not always exactly
so, and asked if it is nevertheless true that the maximal value of $R_{ij}$ 
is achieved for
some pair of vertices $\setof{i,j}$ whose distance apart in the graph is 
maximal.  Here is a counter-example.

\begin{eg}
Let $k=2$, so that the block design is the same as its concurrence graph.
Take $v=10$ and $b=14$.  The graph consists of a cube, with two extra vertices
$1$ and $2$ attached as leaves to vertex $3$.  The vertex antipodal to~$3$ in 
the cube is labelled~$4$.  It is straightforward to check (either using an 
electrical network, or by using the fact that the cube is distance-regular) 
that the effective resistance between a pair of cube vertices is $7/12$, 
$3/4$ and $5/6$ for vertices at distances $1$, $2$ and $3$. 
Hence $R_{1j}\leq 11/6$ for all cube vertices~$j$, while $R_{12}=2$.
On the other hand, the distance between vertices $1$ and $2$ is only~$2$, 
while that between either of them and vertex~$4$ is $4$.
\end{eg}

There are some `nice' graphs where pairwise variance does indeed
increase with distance.  The following result is proved in
\cite{rab:dgh}.  Biggs gave the equivalent result for effective
resistances in \cite{biggselec}.

\begin{thm}
Suppose that a block design has just two distinct concurrences, 
and that the pairs of
vertices corresponding to the larger concurrence form the edges of a
distance-regular graph $H$.  Then pairwise variance increases with
distance in $H$.
\end{thm}

\section{Graph concepts linked to E-optimality}

\subsection{Measures of bottlenecks}
A `good' graph (for use as a network) is  one without bottlenecks: 
any set of vertices should have many edges joining it to its complement.
So, for any subset $S$ of vertices, we let $\partial(S)$
(the \emph{boundary} of $S$) be the set of edges which have one vertex in
$S$ and the other in its complement, and then define the
\emph{isoperimetric number} $\iota(G)$ by
\[
\iota(G) = \min\left\{\frac{|\partial S|}{|S|}: 
S\subseteq V(G),\ 0<|S|\le \frac{v}{2}
\right\}.
\]

The next result shows that the isoperimetric
number is related to the E-criterion. It is useful not 
so much for identifying the E-optimal designs as for easily showing that large
classes of designs cannot be E-optimal: any design whose concurrence graph has
low isoperimetric number performs poorly on the E-criterion.

\begin{cutset}
Let $G$ have an edge-cutset of size~$c$ whose removal separates the graph
into 
parts $S$ and $G\setminus S$ with $m$ and $n$ vertices respectively,
where $0 < m\leq n$. Then
\[
\theta_1\le c\left(\frac{1}{m}+\frac{1}{n}\right) 
\leq \frac{2\left|\partial S\right|}{\left|S\right|}.
\]
\end{cutset}

\begin{pf}
We know that
$\theta_1$ is the minimum of 
$\mathbf{x}^\top \mathbf{L}\mathbf{x}/\mathbf{x}^\top \mathbf{x}$
over real vectors $\mathbf{x}$ with $\sum_i x_i = 0$.
Put 
\[x_i= \left\{\begin{array}{rl}n & \mbox{ if $i\in S$}\\
                          -m & \mbox{otherwise.}
                          \end{array}\right.\]
Then $\mathbf{x}^\top \mathbf{x} = nm(m+n)$ and
\[\mathbf{x}^\top \mathbf{L} \mathbf{x} = 
\sum_{\mathrm{edges\ }ij} (x_i-x_j)^2 = c(m+n)^2.\]
Hence
\[\theta_1 \leq 
\frac{\mathbf{x}^\top \mathbf{L} \mathbf{x}}{\mathbf{x}^\top \mathbf{x}} 
= \frac{c(m+n)^2}{nm(m+n)}
= c\left(\frac{1}{m} + \frac{1}{n}\right) \leq \frac{2c}{m} =
\frac{2\left|\partial S\right|}{\left|S\right|}.
\halmos
\]
\end{pf}

\begin{cor}
\label{thm:iso}
Let $\theta_1$ be the smallest non-trivial eigenvalue of the Laplacian
matrix $\mathbf{L}$ of the connected graph $G$.  Then $\theta_1 \leq 2
\iota(G)$. 
\end{cor}

There is also an upper bound for the isoperimetric number in terms of 
$\theta_1$, which is loosely referred to as a `Cheeger-type inequality';
for details, see the further reading.

We also require a second cutset lemma, phrased in terms of vertex cutsets.

\begin{cutset}
Let $G$ have a vertex-cutset~$C$ of size $c$ whose removal separates the
graph into parts $S$ and $T$ with $m$, $n$ vertices respectively 
(so $nm>0$).  Let $m'$ and $n'$ be the number of edges from vertices in $C$ 
to vertices in $S$, $T$ respectively.  Then
\[
\theta_1 \leq \frac{m'n^2 + n'm^2}{nm(m+n)}.
\]
In particular, if there are no multiple edges at any vertex of $C$ then
$\theta_1\leq c$, with equality if and only if every vertex in $C$ is joined 
to every vertex in $S\cup T$.
\end{cutset}

\begin{pf}
Put
\[x_i= \left\{\begin{array}{rl}n & \mbox{ if $i\in S$}\\
                          -m & \mbox{ if $i\in T$}\\
                           0 & \mbox{ otherwise.}
                          \end{array}\right.\]
Then $\mathbf{x}^\top \mathbf{x} = nm(m+n)$ 
and $\mathbf{x}^\top \mathbf{L}\mathbf{x} = m'n^2+n'm^2$, and so
\[
\theta_1 \leq \frac{m'n^2 + n'm^2}{nm(m+n)}.
\]
If there are no multiple edges at any vertex in $C$ then $m'\leq cm$ and 
$n'\leq cn$ and the result follows.
\halmos
\end{pf}

\subsection{Variance balance}

A block design is \emph{variance-balanced} if all
the concurrences $\lambda_{ij}$ are equal for $i\ne j$.
In such a design, all of the pairwise variances $V_{ij}$ are equal.
Morgan and Srivastav proved the following result in \cite{ms}.

\begin{thm}
If the constant concurrence $\lambda$ of a variance-balanced
design satisfies $(v-1)\lambda = \lfloor(bk/v)\rfloor (k-1)$ then the 
design is E-optimal.
\end{thm}

A block with $k$ different treatments contributes $k(k-1)/2$
edges to the concurrence graph.  Let us define the \emph{defect} of a block 
to be
\[
\frac{k(k-1)}{2} - \mbox{the number of edges it contributes to the graph.}
\]

The following result is proved in \cite{bcc09}.

\begin{thm}
\label{thm:defect}
If $k<v$, then a variance-balanced design with $v$ treatments is E-optimal
if the sum of the block defects is less than $v/2$.
\end{thm}

Table~\ref{tab:binary}(b) shows that the design in Fig.~\ref{fig:binary}(b) 
is variance-balanced.  Block $\Gamma_1$ has defect~$1$, 
and each other block has defect $0$, so the sum of the block defects
is certainly less than $5/2$, and
Theorem~\ref{thm:defect} shows that  the design is E-optimal.  It is rather
counter-intuitive that the non-binary design in Fig.~\ref{fig:binary}(b) can
be better than the design in Fig.~\ref{fig:binary}(a); in fact, in his 
contribution to the discussion of Tocher's paper \cite{tocher}, which 
introduced this design, David Cox said
\begin{quote}
I suspect that \ldots\ balanced ternary designs are of no practical value.
\end{quote}

Computation shows that the design in Fig.~\ref{fig:binary}(a) is 
$\Phi_p$-better than the one in Fig.~\ref{fig:binary}(b) if $p<5.327$.  In 
particular, it is A- and D-better.

\section{Some history}
As we have seen, if the experimental units form a single block and
there are only two treatments then it is 
best for their replications to be as equal as possible.  Statisticians know
this so well that it is hard for us to imagine that more information 
may be obtained, about \emph{all} treatment comparisons, if replications differ
by more than~$1$.

In agriculture, or in any area with qualitative treatments, A-optimality is the
natural criterion.  If treatments are quantities of different substances,
then D-optimality is preferable, as the ranking on this criterion
is invariant to change of measurement units. Thus 
industrial statisticians have tended to prefer  D-optimality,
although E-optimality has become popular among chemical process engineers. 
Perhaps the different camps have not talked to each other as much as they 
should have. 

For most of the 20th century, it was normal practice in field experiments
to have all treatments replicated three or four times.  Where incomplete blocks
were used, they typically had size from $3$ to $20$.  
Yates introduced his square lattice designs with $v=k^2$ in
\cite{fy:lattice}.  He used uniformity data
and two worked examples to show that these designs can give
lower average pairwise variance than a design using a highly
replicated control treatment, but both of his examples were 
equireplicate with $r\in\{3,4\}$.

In the 1930s, 1940s and 1950s, analysis of the data from an experiment
involved inverting the Laplacian matrix without a computer: this is
easy for BIBDs, and only slightly harder if the Laplacian matrix has
only two distinct  non-trivial eigenvalues.  The results in
\cite{kiefer} and \cite{ksh} encouraged the beliefs that the optimal
designs, on all $\Phi_p$-criteria, are as equireplicate as possible, 
with concurrences as equal as possible, and that
the same designs are optimal, or nearly so,  on all of these criteria.

Three short papers in the same journal in 1977--1982 demonstrate
the beliefs at that time.  In \cite{JAJM:rgd}, John and Mitchell did
not even consider designs with unequal replication.  They conjectured
that, if there exist any regular-graph designs for given values of
$v$, $b$ and $k$, then the A- and D-optimal designs are regular-graph
designs.  For the parameter sets which they had examined by computer
search, the same designs were optimal on the A- and D-criteria. In
\cite{BJJAE}, Jones and Eccleston reported the results of various
computer searches for A-optimal designs without the constraint of
equal replication.  For $k=2$ and $b=v\in\{10,11,12\}$ (but not $v=9$) their
A-optimal design is almost a queen-bee design, and their designs are
D-worse than those in \cite{JAJM:rgd}.  The belief in equal
replication was so ingrained that some readers assumed that there was
an error in their program.

John and Williams followed this with the  paper \cite{yyy}
on conjectures for optimal block designs for given values of 
$v$, $b$ and $k$.  Their conjectures included:
\begin{itemize}
\item
the set of regular-graph designs always contains
one that is optimal without this restriction;
\item among regular-graph designs, 
the same designs are optimal on the A- and D-criteria.
\end{itemize}
They endorsed Cox's dismissal  of non-binary designs, strengthening
it to the statement that they ``are inefficient'', and declared that
the three unequally replicated A-optimal designs in \cite{BJJAE} were
``of academic rather than of practical interest''.
These conjectures and opinions seemed quite reasonable to people who 
had been finding good designs for the sizes needed in agricultural experiments.


%
%

At the end of the 20th century, there was an explosion in the number 
of experiments in genomics, using microarrays.  Simplifying the story greatly, 
these are effectively block designs with $k=2$, and biologists wanted 
A-optimal designs, but they did not  know the vocabulary `block' or 
`A-optimal', `graph' or `cycle'.  
Computers were now much more powerful than in 1980, and 
researchers in genomics could simply undertake computer searches without the 
benefit of any statistical theory.
In 2001, Kerr and Churchill \cite{KCh} published the results of a computer
search for A-optimal designs with $k=2$ and $v=b\leq 11$. For
$v\in\{10,11,12\}$, their results were
completely consistent with those in \cite{BJJAE}, which they did not cite.
They called cycles \textit{loop designs}.

Mainstream statisticians began to get involved.
In 2005, Wit, Nobile and Khanin published the paper \cite{Wit} giving the
results of a computer search for A- and D-optimal designs with $k=2$ and 
$v=b$.  The results are shown in Fig.~\ref{fig:wit}.
The A-optimal designs differ from the D-optimal designs
when $v\geq 9$, but are consistent with those found in \cite{KCh}.

\begin{figure}
\begin{center}
\begin{tabular}{|c|c|c|}
\hline
\setlength{\unitlength}{0.4pt}
\begin{picture}(200,200)(-100,-100)
\put(85.1,0){\blobbb}
\put(42.6,73.7){\blobbb}
\put(-42.6,73.7){\blobbb}
\put(-85.1,0){\blobbb}
\put(-42.6,-73.7){\blobbb}
\put(42.6,-73.7){\blobbb}
\qbezier(85.1,0)(63.85,36.85)(42.6,73.7)
\qbezier(42.6,73.7)(0,73.7)(-42.6,73.7)
\qbezier(-85.1,0)(-63.85,36.85)(-42.6,73.7)
\qbezier(85.1,0)(63.85,-36.85)(42.6,-73.7)
\qbezier(42.6,-73.7)(0,-73.7)(-42.6,-73.7)
\qbezier(-85.1,0)(-63.85,-36.85)(-42.6,-73.7)
\end{picture}
& 
\setlength{\unitlength}{0.4pt}
\begin{picture}(200,200)(-100,-100)
\put(85.1,0){\blobbb}
\put(53.1,66.5){\blobbb}
\put(-18.9,83){\blobbb}
\put(-76.7,36.9){\blobbb}
\put(-76.7,-36.9){\blobbb}
\put(-18.9,-83){\blobbb}
\put(53.1,-66.5){\blobbb}
\qbezier(85.1,0)(69.1,33.25)(53.1,66.5)
\qbezier(53.1,66.5)(17.1,74.75)(-18.9,83)
\qbezier(-18.9,83)(-47.8,59.95)(-76.7,36.9)
\qbezier(-76.7,36.9)(-76.7,0)(-76.7,-36.9)
\qbezier(-18.9,-83)(-47.8,-59.95)(-76.7,-36.9)
\qbezier(53.1,-66.5)(17.1,-74.75)(-18.9,-83)
\qbezier(85.1,0)(69.1,-33.25)(53.1,-66.5)
\end{picture}
&
\setlength{\unitlength}{0.4pt}
\begin{picture}(200,200)(-100,-100)
\put(85.1,0){\blobbb}
\put(60.2,60.2){\blobbb}
\put(0,85.1){\blobbb}
\put(-60.2,60.2){\blobbb}
\put(-85.1,0){\blobbb}
\put(-60.2,-60.2){\blobbb}
\put(0,-85.1){\blobbb}
\put(60.2,-60.2){\blobbb}
\qbezier(85.1,0)(72.75,30.1)(60.2,60.2)
\qbezier(60.2,60.2)(30.1,72.75)(0,85.1)
\qbezier(-85.1,0)(-72.75,30.1)(-60.2,60.2)
\qbezier(-60.2,60.2)(-30.1,72.75)(0,85.1)
\qbezier(85.1,0)(72.75,-30.1)(60.2,-60.2)
\qbezier(60.2,-60.2)(30.1,-72.75)(0,-85.1)
\qbezier(-85.1,0)(-72.75,-30.1)(-60.2,-60.2)
\qbezier(-60.2,-60.2)(-30.1,-72.75)(0,-85.1)
\end{picture}
\\
\hline
\setlength{\unitlength}{0.4pt}
\begin{picture}(200,200)(-100,-100)
\put(85.1,0){\blobbb}
\put(65.2,54.7){\blobbb}
\put(14.8,83.8){\blobbb}
\put(-42.6,73.7){\blobbb}
\put(-80,29.1){\blobbb}
\put(-80,-29.1){\blobbb}
\put(-42.6,-73.7){\blobbb}
\put(14.8,-83.8){\blobbb}
\put(65.2,-54.7){\blobbb}
\qbezier(85.1,0)(75.15,27.35)(65.2,54.7)
\qbezier(65.2,54.7)(40,69.25)(14.8,83.8)
\qbezier(14.8,83.8)(-13.9,78.75)(-42.6,73.7)
\qbezier(-42.6,73.7)(-62.3,51.4)(-80,29.1)
\qbezier(-80,29.1)(-80,0)(-80,-29.1)
\qbezier(85.1,0)(75.15,-27.35)(65.2,-54.7)
\qbezier(65.2,-54.7)(40,-69.25)(14.8,-83.8)
\qbezier(14.8,-83.8)(-13.9,-78.75)(-42.6,-73.7)
\qbezier(-42.6,-73.7)(-62.3,-51.4)(-80,-29.1)
\end{picture}
&
\setlength{\unitlength}{0.4pt}
\begin{picture}(200,200)(-100,-31.2)
\put(50,0){\blobbb}
\put(80.9,95.1){\blobbb}
\put(0,153.9){\blobbb}
\put(-80.9,95.1){\blobbb}
\put(-50,0){\blobbb}
\put(0,-16.3){\blobbb}
\put(80.9,42.5){\blobbb}
\put(50,137.6){\blobbb}
\put(-50,137.6){\blobbb}
\put(-80.9,42.5){\blobbb}
\qbezier(50,0)(65.45,21.25)(80.9,42.5)
\qbezier(-50,0)(-65.45,21.25)(-80.9,42.5)
\qbezier(80.9,95.1)(65.45,118.35)(50,137.6)
\qbezier(-80.9,95.1)(-65.45,118.35)(-50,137.6)
\qbezier(0,153.9)(-25,145.75)(-50,137.6)
\qbezier(0,153.9)(25,145.75)(50,137.6)
\qbezier(-80.9,95.1)(-80.9,68.8)(-80.9,42.5)
\qbezier(80.9,95.1)(80.9,68.8)(80.9,42.5)
\qbezier(0,-16.3)(25,-8.15)(50,0)
\qbezier(0,-16.3)(-25,-8.15)(-50,0)
\end{picture}
&
\setlength{\unitlength}{0.4pt}
\begin{picture}(200,200)(-100,-100)
\put(85.1,0){\blobbb}
\put(71.6,46){\blobbb}
\put(35.4,77.4){\blobbb}
\put(-12.1,84.2){\blobbb}
\put(-55.7,64.3){\blobbb}
\put(-81.7,24){\blobbb}
\put(71.6,-46){\blobbb}
\put(35.4,-77.4){\blobbb}
\put(-12.1,-84.2){\blobbb}
\put(-55.7,-64.3){\blobbb}
\put(-81.7,-24){\blobbb}
\qbezier(85.1,0)(78.35,23)(71.6,46)
\qbezier(71.6,46)(53.5,61.7)(35.4,77.4)
\qbezier(35.4,77.4)(11.65,80.8)(-12.1,84.2)
\qbezier(-12.1,84.2)(-33.9,74.25)(-55.7,64.3)
\qbezier(-55.7,64.3)(-68.7,44.15)(-81.7,24)
\qbezier(-81.7,24)(-81.7,0)(-81.7,-24)
\qbezier(85.1,0)(78.35,-23)(71.6,-46)
\qbezier(71.6,-46)(53.5,-61.7)(35.4,-77.4)
\qbezier(35.4,-77.4)(11.65,-80.8)(-12.1,-84.2)
\qbezier(-12.1,-84.2)(-33.9,-74.25)(-55.7,-64.3)
\qbezier(-55.7,-64.3)(-68.7,-44.15)(-81.7,-24)
\end{picture}
\\
\hline
\multicolumn{3}{c}{\rule[-20pt]{0pt}{40pt}(a) D-optimal designs}\\
\hline
\setlength{\unitlength}{0.4pt}
\begin{picture}(200,200)(-100,-100)
\put(85.1,0){\blobbb}
\put(42.6,73.7){\blobbb}
\put(-42.6,73.7){\blobbb}
\put(-85.1,0){\blobbb}
\put(-42.6,-73.7){\blobbb}
\put(42.6,-73.7){\blobbb}
\qbezier(85.1,0)(63.85,36.85)(42.6,73.7)
\qbezier(42.6,73.7)(0,73.7)(-42.6,73.7)
\qbezier(-85.1,0)(-63.85,36.85)(-42.6,73.7)
\qbezier(85.1,0)(63.85,-36.85)(42.6,-73.7)
\qbezier(42.6,-73.7)(0,-73.7)(-42.6,-73.7)
\qbezier(-85.1,0)(-63.85,-36.85)(-42.6,-73.7)
\end{picture}
&
\setlength{\unitlength}{0.4pt}
\begin{picture}(200,200)(-100,-100)
\put(85.1,0){\blobbb}
\put(53.1,66.5){\blobbb}
\put(-18.9,83){\blobbb}
\put(-76.7,36.9){\blobbb}
\put(-76.7,-36.9){\blobbb}
\put(-18.9,-83){\blobbb}
\put(53.1,-66.5){\blobbb}
\qbezier(85.1,0)(69.1,33.25)(53.1,66.5)
\qbezier(53.1,66.5)(17.1,74.75)(-18.9,83)
\qbezier(-18.9,83)(-47.8,59.95)(-76.7,36.9)
\qbezier(-76.7,36.9)(-76.7,0)(-76.7,-36.9)
\qbezier(-18.9,-83)(-47.8,-59.95)(-76.7,-36.9)
\qbezier(53.1,-66.5)(17.1,-74.75)(-18.9,-83)
\qbezier(85.1,0)(69.1,-33.25)(53.1,-66.5)
\end{picture}
&
\setlength{\unitlength}{0.4pt}
\begin{picture}(200,200)(-100,-100)
\put(85.1,0){\blobbb}
\put(60.2,60.2){\blobbb}
\put(0,85.1){\blobbb}
\put(-60.2,60.2){\blobbb}
\put(-85.1,0){\blobbb}
\put(-60.2,-60.2){\blobbb}
\put(0,-85.1){\blobbb}
\put(60.2,-60.2){\blobbb}
\qbezier(85.1,0)(72.75,30.1)(60.2,60.2)
\qbezier(60.2,60.2)(30.1,72.75)(0,85.1)
\qbezier(-85.1,0)(-72.75,30.1)(-60.2,60.2)
\qbezier(-60.2,60.2)(-30.1,72.75)(0,85.1)
\qbezier(85.1,0)(72.75,-30.1)(60.2,-60.2)
\qbezier(60.2,-60.2)(30.1,-72.75)(0,-85.1)
\qbezier(-85.1,0)(-72.75,-30.1)(-60.2,-60.2)
\qbezier(-60.2,-60.2)(-30.1,-72.75)(0,-85.1)
\end{picture}
\\
\hline
\setlength{\unitlength}{0.4pt}
\begin{picture}(130,130)(-65,-65)
\put(0,0){\blobbb}
\put(-50,0){\blobbb}
\put(50,0){\blobbb}
\put(-50,-50){\blobbb}
\put(0,-50){\blobbb}
\put(0,50){\blobbb}
\put(43.3,25){\blobbb}
\put(43.3,-25){\blobbb}
\put(25,43.3){\blobbb}
\qbezier(0,0)(-25,0)(-50,0)
\qbezier(0,0)(0,-25)(0,-50)
\qbezier(0,0)(0,25)(0,50)
\qbezier(0,0)(25,0)(50,0)
\qbezier(-50,-50)(-50,-25)(-50,0)
\qbezier(-50,-50)(-25,-50)(0,-50)
\qbezier(0,0)(21.65,-12.5)(43.3,-25)
\qbezier(0,0)(21.65,12.5)(43.3,25)
\qbezier(0,0)(12.5,21.65)(25,43.3)
\end{picture}
& 
\setlength{\unitlength}{0.4pt}
\begin{picture}(130,130)(-65,-65)
\put(0,0){\blobbb}
\put(-50,0){\blobbb}
\put(50,0){\blobbb}
\put(-50,-50){\blobbb}
\put(0,-50){\blobbb}
\put(0,50){\blobbb}
\put(43.3,25){\blobbb}
\put(43.3,-25){\blobbb}
\put(25,43.3){\blobbb}
\put(-25,43.3){\blobbb}
\qbezier(0,0)(-25,0)(-50,0)
\qbezier(0,0)(0,-25)(0,-50)
\qbezier(0,0)(0,25)(0,50)
\qbezier(0,0)(25,0)(50,0)
\qbezier(-50,-50)(-50,-25)(-50,0)
\qbezier(-50,-50)(-25,-50)(0,-50)
\qbezier(0,0)(21.65,-12.5)(43.3,-25)
\qbezier(0,0)(21.65,12.5)(43.3,25)
\qbezier(0,0)(12.5,21.65)(25,43.3)
\qbezier(0,0)(-12.5,21.65)(-25,43.3)
\end{picture}
& 
\setlength{\unitlength}{0.4pt}
\begin{picture}(130,130)(-65,-65)
\put(0,0){\blobbb}
\put(-50,0){\blobbb}
\put(50,0){\blobbb}
\put(-50,-50){\blobbb}
\put(0,-50){\blobbb}
\put(0,50){\blobbb}
\put(43.3,25){\blobbb}
\put(43.3,-25){\blobbb}
\put(25,43.3){\blobbb}
\put(-25,43.3){\blobbb}
\put(-43.3,25){\blobbb}
\qbezier(0,0)(-25,0)(-50,0)
\qbezier(0,0)(0,-25)(0,-50)
\qbezier(0,0)(0,25)(0,50)
\qbezier(0,0)(25,0)(50,0)
\qbezier(-50,-50)(-50,-25)(-50,0)
\qbezier(-50,-50)(-25,-50)(0,-50)
\qbezier(0,0)(21.65,-12.5)(43.3,-25)
\qbezier(0,0)(21.65,12.5)(43.3,25)
\qbezier(0,0)(12.5,21.65)(25,43.3)
\qbezier(0,0)(-12.5,21.65)(-25,43.3)
\qbezier(0,0)(-21.65,12.5)(-43.3,25)
\end{picture}
\\
\hline
\multicolumn{3}{c}{\rule[-0pt]{0pt}{20pt}(b) A-optimal designs}\\
\end{tabular}
\end{center}
\caption{D-and A-optimal designs with $k=2$ and $6\leq v=b\leq 11$}
\label{fig:wit}
\end{figure}

What is going on here?  Why are the designs so different when $v\geq 9$?  Why 
is there such a sudden, large change in the A-optimal designs?  We explain this
in the next section.

\section{Block size two}
\subsection{Least replication}
\label{sec:least}

\begin{figure}
\begin{center}
\begin{tabular}{c@{\qquad\qquad}c}
\setlength{\unitlength}{0.8pt}
\begin{picture}(100,100)
\put(50,50){\blobbb}
\put(50,100){\blobbb}
\put(50,0){\blobbb}
\put(0,50){\blobbb}
\put(100,50){\blobbb}
\put(83.3,83.3){\blobbb}
\put(16.7,83.3){\blobbb}
\put(83.3,16.7){\blobbb}
\put(16.7,16.7){\blobbb}
\put(50,0){\line(-2,1){33.3}}
\put(0,50){\line(1,0){50}}
\put(0,50){\line(1,-2){16.7}}
\put(100,50){\line(-1,-2){16.7}}
\put(100,50){\line(-1,2){16.7}}
\put(50,100){\line(-2,-1){33.3}}
\put(16.7,83.3){\line(1,-1){67}}
\end{picture}
&
\setlength{\unitlength}{0.8pt}
\begin{picture}(100,100)
\put(50,50){\blobbb}
\put(50,100){\blobbb}
\put(50,0){\blobbb}
\put(0,50){\blobbb}
\put(100,50){\blobbb}
\put(85,85){\blobbb}
\put(15,85){\blobbb}
\put(85,15){\blobbb}
\put(15,15){\blobbb}
\put(0,50){\line(1,0){100}}
\put(50,0){\line(0,1){100}}
\put(15,15){\line(1,1){70}}
\put(15,85){\line(1,-1){70}}
\end{picture}
\end{tabular}
\end{center}
\caption{Two trees with $v=9$, $b=8$ and $k=2$}
\label{fig:tree}
\end{figure}
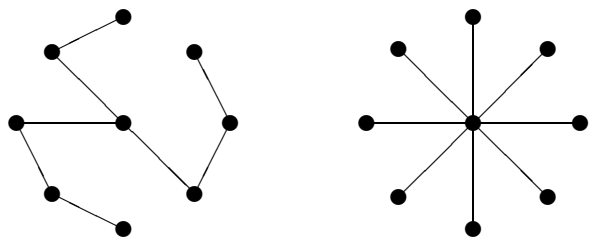

If $k=2$, then the design is the same as its concurrence graph,
and connectivity requires that $b\geq v-1$.
If $b=v-1$, then all connected designs are trees, such as those in 
Fig.~\ref{fig:tree}.
Theorem~\ref{thm:span} shows that the D-criterion does not 
differentiate between them.

In a tree, the effective resistance $R_{ij}$ is just the length of the unique
path between vertices $i$ and $j$.  Theorems~\ref{thm:var} and~\ref{thm:resist}
show that the only A-optimal designs are the stars, such as the graph on the
right of Fig.~\ref{fig:tree}.

In a star with $v$~vertices, the contrast between any two leaves is an 
eigenvector of the Laplacian matrix~$\mathbf{L}$ with eigenvalue~$1$, 
while the 
contrast between the central vertex and all the other vertices is an 
eigenvector with eigenvalue~$v$. If $v\geq5$ and $G$ is not a star then there
is an edge whose removal splits the graph into two components of sizes at least
$2$ and $3$.  Cutset Lemma~1 then shows that $\theta_1 \leq 5/6 <1$.  The
only other tree which is not a star is the path of length $3$, for which 
direct calculation shows that $\theta_1 = 2 - \sqrt{2} <1$.  Hence the 
E-optimal designs are also the stars.

\subsection{One fewer treatment}
\label{sec:loop}

\begin{figure}[htbp]
\begin{center}
\begin{tabular}{@{}c@{\qquad\qquad}c@{}}
\setlength{\unitlength}{8mm}
\begin{picture}(8,7.5)(-2,-2)
\put(2,0){\blobb}
\put(-2,0){\blobb}
\put(1,1.67){\blobb}
\put(3,1.67){\blobb}
\put(2,3.33){\blobb}
\put(4,3.33){\blobb}
\put(4,5){\blobb}
\put(6,3.33){\blobb}
\put(-1,1.67){\blobb}
\put(-2,3.33){\blobb}
\put(1,-1.67){\blobb}
\put(-1,-1.67){\blobb}
\put(0.7,1.37){\makebox(0,0){{{$6$}}}}
\put(-0.7,1.37){\makebox(0,0){{{$1$}}}}
\put(-1.6,0){\makebox(0,0){{$2$}}}
\put(-0.7,-1.37){\makebox(0,0){{$3$}}}
\put(0.7,-1.37){\makebox(0,0){{$4$}}}
\put(1.6,0){\makebox(0,0){{$5$}}}
\put(-1.4,3.33){\makebox(0,0){{$12$}}}
\put(3.4,1.67){\makebox(0,0){{{$7$}}}}
\put(1.6,3.33){\makebox(0,0){{{$8$}}}}
\put(4,2.83){\makebox(0,0){{$9$}}}
\put(3.4,5){\makebox(0,0){{$10$}}}
\put(6,2.83){\makebox(0,0){{$11$}}}
\put(1,1.67){\line(-1,0){2}}
\put(1,-1.67){\line(-1,0){2}}
\put(3,1.67){\line(-1,0){2}}
\put(4,3.33){\line(-1,0){2}}
\put(6,3.33){\line(-1,0){2}}
\put(4,3.33){\line(0,1){1.63}}
\put(-1,1.67){\line(-3,5){1}}
\put(-1,-1.67){\line(-3,5){1}}
\put(2,0){\line(-3,5){1}}
\put(1,1.67){\line(3,5){1}}
\put(1,-1.67){\line(3,5){1}}
\put(-2,0){\line(3,5){1}}
\end{picture}
&
\setlength{\unitlength}{8mm}
\begin{picture}(6,7.5)(-2,-2)
\put(2,0){\blobb}
\put(-2,0){\blobb}
\put(1,1.67){\blobb}
\put(3,1.67){\blobb}
\put(2,3.33){\blobb}
\put(4,3.33){\blobb}
\put(2,5){\blobb}
\put(0,3.33){\blobb}
\put(-1,1.67){\blobb}
\put(-2,3.33){\blobb}
\put(1,-1.67){\blobb}
\put(-1,-1.67){\blobb}
\put(0.7,1.37){\makebox(0,0){{{$6$}}}}
\put(-0.7,1.37){\makebox(0,0){{{$1$}}}}
\put(-1.6,0){\makebox(0,0){{$2$}}}
\put(-0.7,-1.37){\makebox(0,0){{$3$}}}
\put(0.7,-1.37){\makebox(0,0){{$4$}}}
\put(1.6,0){\makebox(0,0){{$5$}}}
\put(-1.4,3.33){\makebox(0,0){{$12$}}}
\put(3.4,1.67){\makebox(0,0){{{$7$}}}}
\put(0.4,3.33){\makebox(0,0){{{$8$}}}}
\put(2,2.83){\makebox(0,0){{$9$}}}
\put(1.4,5){\makebox(0,0){{$10$}}}
\put(4,2.83){\makebox(0,0){{$11$}}}
\put(1,1.67){\line(-1,0){2}}
\put(1,-1.67){\line(-1,0){2}}
\put(3,1.67){\line(-1,0){2}}
\put(4,3.33){\line(-1,0){2}}
\put(2,3.33){\line(0,1){1.63}}
\put(-1,1.67){\line(-3,5){1}}
\put(-1,-1.67){\line(-3,5){1}}
\put(2,0){\line(-3,5){1}}
\put(1,1.67){\line(3,5){1}}
\put(1,1.67){\line(-3,5){1}}
\put(1,-1.67){\line(3,5){1}}
\put(-2,0){\line(3,5){1}}
\end{picture}
\\ \\
(a) & (b)
\end{tabular}
\end{center}
\caption{Two unicyclic graphs with $b=v=12$ and $s=6$}
\label{fig:move}
\end{figure}
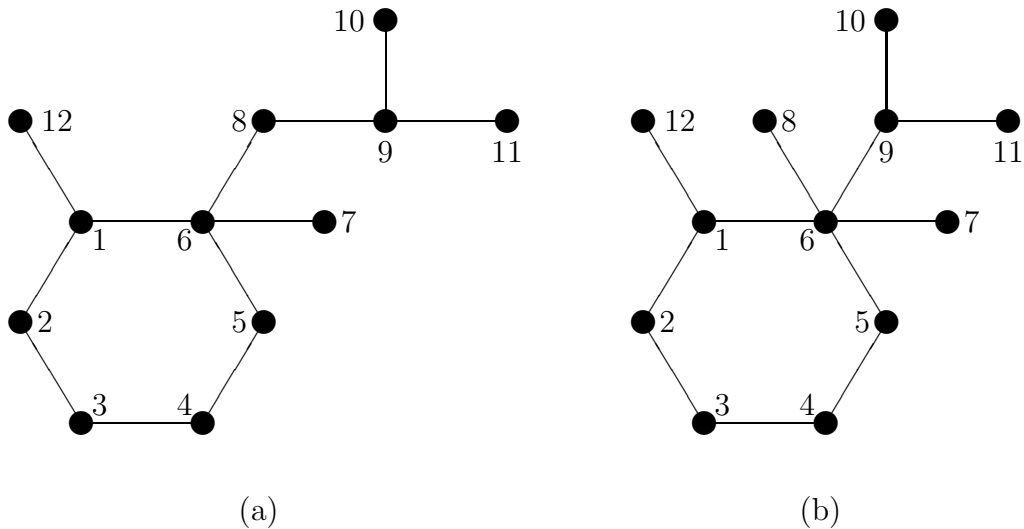

If $b=v$ and $k=2$, then the concurrence graph $G$ contains a single cycle:
such graphs are called \emph{unicyclic}. Let $s$ be the length of the cycle.
All the remaining vertices are in trees attached to various vertices of 
the cycle.
Fig.~\ref{fig:move} shows two unicyclic graphs with $v=12$ and $s=6$.

As we remarked in Section~\ref{sec:span}, the number of spanning trees in a
unicyclic graph is equal to the length of the cycle. Hence, 
Theorem~\ref{thm:span} gives the following result.

\begin{thm}
If $k=2$ and $b=v\geq3$, then the D-optimal designs are precisely the cycles.
\end{thm}

For A-optimality, we first show that no graph like the one in 
Fig.~\ref{fig:move}(a) can be optimal.  If vertex~$12$ is moved so that
it is joined to vertex~$6$, instead of vertex~$1$, then the sum of the 
variances $V_{i,12}$ for $i$ in the cycle is unchanged and the variances 
$V_{i,12}$ for the remaining vertices $i$ are all decreased.  This argument 
shows that all the trees must be attached to the same vertex of the cycle.

Now consider the tree on vertices $6$, $8$, $9$, $10$ and $11$ in 
Fig.~\ref{fig:move}(a). If the two edges incident with vertex~$8$ are
modified to those in Fig.~\ref{fig:move}(b), then the set of
variances between these five vertices are unchanged, as are all others
involving vertex~$8$, but those between vertices $9$, $10$, $11$ and vertices
outside this tree are all decreased.  This argument shows that, for any given
length~$s$ of the cycle, the only candidate for an A-optimal design has $v-s$ 
leaves  attached to a single vertex of the cycle.

The effective resistance between a pair of vertices at distance $d$ in a cycle
of length~$s$ is $d(s-d)/s$, while that between a leaf and the cycle vertex to
which it is attached is $1$.  Hence a short calculation shows that the sum of
the pairwise effective resistances is equal to $g(s)/12$,
where 
\[g(s) = -s^3 +2vs^2 + 13s -12sv +12v^2 -14v.\] 
Now $\bar V/\sigma^2 = g(s)/[3v(v-1)]$ and
we seek the minimum of $g(s)$ for integers $s$ in the interval $[2,v]$.

\begin{figure}
\begin{center}
\setlength{\unitlength}{1.9cm}
\begin{picture}(7,3.6)(-0.6,1.8)
\put(0,2){\line(1,0){5.5}}
\put(0,2){\line(0,1){3.2}}
\xxxlabel{1}{3}
\xxxlabel{2}{6}
\xxxlabel{3}{9}
\xxxlabel{4}{12}
\xxxlabel{5}{15}
\yyylabel{2}{2}
\yyylabel{3}{3}
\yyylabel{4}{4}
\yyylabel{5}{5}
\put(5.8,2){\makebox(0,0){$s$}}
\put(0,5.35){\makebox(0,0){$\bar V/\sigma^2$}}
\put(0.67,3){\makebox(0,0){$\boldsymbol{\star}$}}
\put(1,2.8){\makebox(0,0){$\boldsymbol{\star}$}}
\put(1.33,2.67){\makebox(0,0){$\boldsymbol{\star}$}}
\put(1.67,2.53){\makebox(0,0){$\boldsymbol{\star}$}}
\put(2,2.33){\makebox(0,0){$\boldsymbol{\star}$}}
\put(0.67,3.14){\makebox(0,0){$\boldsymbol{\odot}$}}
\put(1,2.98){\makebox(0,0){$\boldsymbol{\odot}$}}
\put(1.33,2.90){\makebox(0,0){$\boldsymbol{\odot}$}}
\put(1.67,2.86){\makebox(0,0){$\boldsymbol{\odot}$}}
\put(2,2.79){\makebox(0,0){$\boldsymbol{\odot}$}}
\put(2.33,2.67){\makebox(0,0){$\boldsymbol{\odot}$}}
\put(0.67,3.25){\makebox(0,0){$\boldsymbol{\ast}$}}
\put(1,3.12){\makebox(0,0){$\boldsymbol{\ast}$}}
\put(1.33,3.07){\makebox(0,0){$\boldsymbol{\ast}$}}
\put(1.67,3.07){\makebox(0,0){$\boldsymbol{\ast}$}}
\put(2,3.083){\makebox(0,0){$\boldsymbol{\ast}$}}
\put(2.33,3.07){\makebox(0,0){$\boldsymbol{\ast}$}}
\put(2.67,3){\makebox(0,0){$\boldsymbol{\ast}$}}
\put(0.67,3.33){\makebox(0,0){$\boldsymbol{\diamond}$}}
\put(1,3.22){\makebox(0,0){$\boldsymbol{\diamond}$}}
\put(1.33,3.19){\makebox(0,0){$\boldsymbol{\diamond}$}}
\put(1.67,3.22){\makebox(0,0){$\boldsymbol{\diamond}$}}
\put(2,3.28){\makebox(0,0){$\boldsymbol{\diamond}$}}
\put(2.33,3.33){\makebox(0,0){$\boldsymbol{\diamond}$}}
\put(2.67,3.36){\makebox(0,0){$\boldsymbol{\diamond}$}}
\put(3,3.33){\makebox(0,0){$\boldsymbol{\diamond}$}}
\put(0.67,3.4){\makebox(0,0){$\boldsymbol{\bullet}$}}
\put(1,3.3){\makebox(0,0){$\boldsymbol{\bullet}$}}
\put(1.33,3.29){\makebox(0,0){$\boldsymbol{\bullet}$}}
\put(1.67,3.33){\makebox(0,0){$\boldsymbol{\bullet}$}}
\put(2,3.41){\makebox(0,0){$\boldsymbol{\bullet}$}}
\put(2.33,3.51){\makebox(0,0){$\boldsymbol{\bullet}$}}
\put(2.67,3.6){\makebox(0,0){$\boldsymbol{\bullet}$}}
\put(3,3.66){\makebox(0,0){$\boldsymbol{\bullet}$}}
\put(3.33,3.67){\makebox(0,0){$\boldsymbol{\bullet}$}}
\put(0.67,3.45){\makebox(0,0){$\boldsymbol{\circ}$}}
\put(1,3.37){\makebox(0,0){$\boldsymbol{\circ}$}}
\put(1.33,3.36){\makebox(0,0){$\boldsymbol{\circ}$}}
\put(1.67,3.42){\makebox(0,0){$\boldsymbol{\circ}$}}
\put(2,3.52){\makebox(0,0){$\boldsymbol{\circ}$}}
\put(2.33,3.64){\makebox(0,0){$\boldsymbol{\circ}$}}
\put(2.67,3.76){\makebox(0,0){$\boldsymbol{\circ}$}}
\put(3,3.88){\makebox(0,0){$\boldsymbol{\circ}$}}
\put(3.33,3.96){\makebox(0,0){$\boldsymbol{\circ}$}}
\put(3.67,4){\makebox(0,0){$\boldsymbol{\circ}$}}
\put(0.67,3.5){\makebox(0,0){$\boldsymbol{\times}$}}
\put(1,3.424){\makebox(0,0){$\boldsymbol{\times}$}}
\put(1.33,3.424){\makebox(0,0){$\boldsymbol{\times}$}}
\put(1.67,3.48){\makebox(0,0){$\boldsymbol{\times}$}}
\put(2,3.59){\makebox(0,0){$\boldsymbol{\times}$}}
\put(2.33,3.73){\makebox(0,0){$\boldsymbol{\times}$}}
\put(2.67,3.88){\makebox(0,0){$\boldsymbol{\times}$}}
\put(3,4.03){\makebox(0,0){$\boldsymbol{\times}$}}
\put(3.33,4.17){\makebox(0,0){$\boldsymbol{\times}$}}
\put(3.67,4.27){\makebox(0,0){$\boldsymbol{\times}$}}
\put(4,4.33){\makebox(0,0){$\boldsymbol{\times}$}}
\put(0.67,3.54){\makebox(0,0){$\boldsymbol{+}$}}
\put(1,3.419){\makebox(0,0){$\boldsymbol{+}$}}
\put(1.33,3.47){\makebox(0,0){$\boldsymbol{+}$}}
\put(1.67,3.54){\makebox(0,0){$\boldsymbol{+}$}}
\put(2,3.65){\makebox(0,0){$\boldsymbol{+}$}}
\put(2.33,3.79){\makebox(0,0){$\boldsymbol{+}$}}
\put(2.67,3.96){\makebox(0,0){$\boldsymbol{+}$}}
\put(3,4.14){\makebox(0,0){$\boldsymbol{+}$}}
\put(3.33,4.31){\makebox(0,0){$\boldsymbol{+}$}}
\put(3.67,4.46){\makebox(0,0){$\boldsymbol{+}$}}
\put(4,4.59){\makebox(0,0){$\boldsymbol{+}$}}
\put(4.33,4.67){\makebox(0,0){$\boldsymbol{+}$}}
\put(4.83,3.73){\makebox{$\begin{array}{cc}
\boldsymbol{+} & v=13\\
\boldsymbol{\times} & v=12\\
\boldsymbol{\circ} & v=11\\
\boldsymbol{\bullet} & v=10\\
\boldsymbol{\diamond} & v=9\\
\boldsymbol{\ast} & v=8\\
\boldsymbol{\odot} & v=7\\
\boldsymbol{\star} & v=6\\
\end{array}$}}
\end{picture}
\end{center}
\caption{Average pairwise variance, in a unicylic graph with $v$~vertices, 
as a function of the length~$s$ of the cycle}
\label{fig:cubic}
\end{figure}
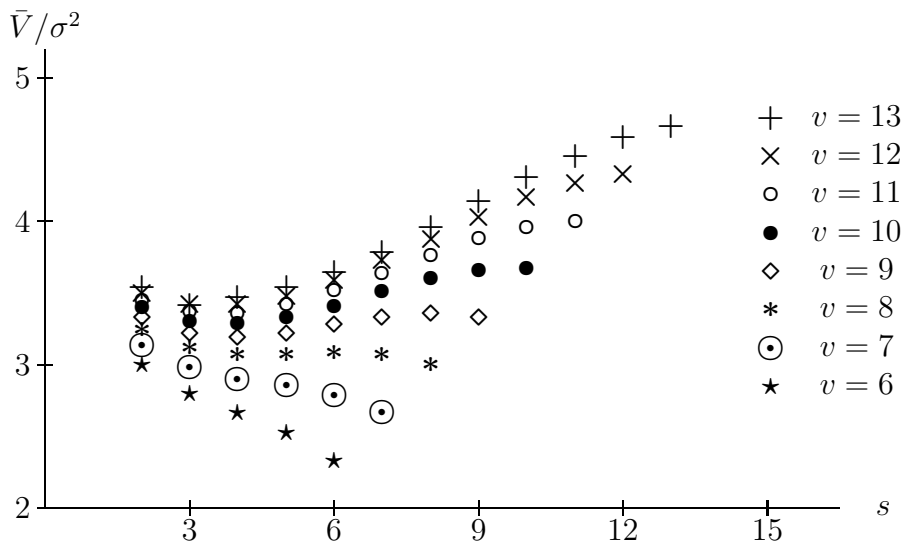

Fig.~\ref{fig:cubic} plots $g(s)/[3v(v-1)]$ for $s$ in $[2,v]$ and 
$6\leq v\leq 13$.
When $v\leq 7$, the function $g$ is monotonic decreasing, so it attains its
minimum on $[2,v]$ at $s=v$. For all larger values of $v$, the 
function~$g$ has a local minimum in the interval $[3,5]$: 
when $v\geq 9$, the value at this local minimum is
less than $g(v)$.  This change from the upper end of the interval to the local 
minimum explains the sudden change in the A-optimal designs.  Detailed 
examination of the local minimum gives the following result.

\begin{thm}
\label{thm:aopt}
If $k=2$ and $b=v\geq3$ then the A-optimal designs are:
\begin{itemize}
\item  a cycle, if $v\leq 8$;
\item a square with $v-4$ leaves attached to one vertex, if $9\leq v \leq 11$;
\item a triangle with $v-4$ leaves attached to one vertex, if $ v \geq 13$;
\item either of the last two, if $v=12$.
\end{itemize}
\end{thm}

What about E-optimality?  The smallest eigenvalue of the Laplacian matrix of
the triangle with one or more leaves attached to one vertex is 1, as is that of
the digon with two or more leaves attached to one vertex.  We now show that
almost all other unicyclic graphs have at least one non-trivial eigenvalue
smaller than this.

Suppose that vertex~$i$ in the cycle has a non-empty tree 
attached to it, so that $\{i\}$ is a vertex cutset.  
If $s\geq 3$ then there are no double edges, so 
Cutset Lemma~2  shows that $\theta_1 <1$ unless all vertices are joined 
to~$i$, in which case $s=3$.
If $s=2$ and there are trees attached to both vertices of the digon, then
 applying Cutset Lemma~2 at each of these vertices shows that
$\theta_1<1$ unless $v=4$ and there is one leaf at each vertex of the digon: 
for this graph, $\theta_1 = 2-\sqrt{5} <1$.  A digon with leaves attached to 
one vertex is just a star with one edge doubled.

The cycle of size~$v$ is a cyclic design.  The smallest eigenvalue of its
Laplacian matrix is $2(1 - \cos(2\pi/v))$, which is greater than $1$ when 
$v\leq 5$, is equal to $1$ when $v=6$, and is less than~$1$ when $v\leq 7$.
When $v=3$ it is equal to $3$, which is greater than $3-\sqrt{3}$, which is 
the smallest Laplacian eigenvalue of the digon with one leaf.

Putting all of this together proves the following result.

\begin{thm}
\label{thm:eopt}
If $k=2$ and $b=v\geq3$, then the E-optimal designs are:
\begin{itemize}
\item  a cycle, if $v\leq 5$;
\item a triangle with $v-3$ leaves attached to one vertex, 
or a star with one edge doubled, if $v\geq 7$;
\item either of the last two, if $v=6$.
\end{itemize}
\end{thm}

Thus, for $v\geq 9$, the ranking on the D-criterion is essentially the 
opposite of the ranking on the A- and E-criteria.  The A- and E-optimal 
designs are far from equireplicate.  The change is sudden, not gradual.
These findings were initially quite shocking to statisticians.

\subsection{More blocks}

What happens when $b$ is larger than $v$ but still has the same order of
magnitude?  The following theorems show that the A- and E-optimal
designs are very different from the D-optimal designs when $v$ is
large.   The proofs of Theorems~\ref{thm:thresh} and \ref{thm:eleaf}
are in \cite{rabpaper} and \cite{bcc09} respectively.

\begin{thm}
Let $G$ be the concurrence graph of a connected block design $\Delta$
with $k=2$ and $b\geq v$.  If $\Delta$ is D-optimal then $G$ does not
contain any bridge (an edge cutset of size one): in particular, $G$
contains no leaves.
\end{thm}

\begin{pf}
Suppose that $\{i,j\}$ is an edge-cutset for $G$.  Let $H$ and $K$ be
the parts of $G$ containing $i$ and $j$, respectively.

Since $G$ is not a tree, we may assume that $H$ is not a tree, and so
there is some edge $e$ in $H$ that is not in every spanning tree for
$H$.  Let $n_1$ and $n_2$ be the numbers of spanning trees for $H$ that
include and exclude $e$, respectively, and let $m$ be the number of
spanning trees for $K$.  Every spanning tree for $G$ consists of
spanning trees for $H$ and $K$ together with the edge $\{i,j\}$.
Hence $G$ has $(n_1+n_2)m$ spanning trees.

Let $\ell$ be a vertex on $e$ with $\ell \ne i$.  Form $G'$ from $G$
by removing edge $e$ and inserting the edge $e'$, where $e' =
\{\ell,j\}$. 

Let $T$ and $T'$ be spanning trees for $H$ and $K$ respectively.  If
$T$ does not contain $e$ then $T \cup \{\{i,j\}\} \cup T'$ and $T \cup
\{e'\} \cup T'$ are both spanning trees for $G'$.  If $T$ contains $e$
then $(T\setminus \{e\}) \cup \{\{i,j\}\} \cup \{e'\} \cup T'$ is a
spanning tree for $G'$.  Hence the number of spanning trees for $G'$
is at least $(2n_2 +n_1)m$, which is greater than $(n_1+n_2)m$ because
$n_2\geq 1$. Hence $G$ does not have the maximal number of spanning
trees and so $\Delta$ is not D-optimal.
\halmos
\end{pf}

\begin{thm}
\label{thm:thresh}
Let $c$ be a positive integer.  Then there is a positive integer $v_c$
such that if $b-v=c$ and $v\geq v_c$ then all A-optimal designs with
$k=2$ contain leaves.
\end{thm}

\begin{thm}
\label{thm:eleaf}
If $20 \leq v \leq b\leq 5v/4$ then the concurrence graph for any
E-optimal design with $k=2$ contains leaves.
\end{thm}

Of course, to obtain a BIBD when $k=2$, $b$ needs to be a quadratic
function of $v$.  What happens if $b$ is merely a linear function of
$v$?  In \cite{bcc09} we conjectured that if $b=cv$ for some constant $c$
then there is a threshold result like the one in
Theorem~\ref{thm:thresh}.   However, current work by Robert Johnson
and Mark Walters \cite{JRJMW} suggests something much more
interesting---that there is a constant $C$ with $3<C<4$ such that if 
$b\geq Cv$ and $k=2$ then all A-optimal designs are (nearly) equireplicate, 
and that
random such graphs (in a suitable model) are close to A-optimal with high
probability.
On the other hand if $b\leq Cv$ then a graph consisting of a large almost
equireplicate part
(all degrees 3 and 4 with average degree close to $Cv$) together with a
suitable number of leaves joined to a single vertex is strictly better
than any queen-bee design.

\subsection{A little more history}

The results on D- and A-optimality in Sections~\ref{sec:least} 
and~\ref{sec:loop}
were proved in \cite{rabpaper}, partly 
to put to rest mutterings that the results of \cite{BJJAE,KCh,Wit} found by 
computer search were incorrect.  The results on E-optimality are in 
\cite{bcc09}. 

In spite of the horror with which these results were greeted, it transpired 
that they were not new.  The D- and E-optimal designs for $b=
(v-1)/(k-1)$ were identified in \cite{bap} in 1991.
The A-optimal designs for $k=2$ and $b=v-1$ had been
given in \cite{mandal} in 1991.  Also in 1991, Tjur gave the
A-optimal designs for $k=2$ and $b=v$ in \cite{tj1}: his proof used the Levi 
graph as an electrical network.

A fairly common response to these unexpected results was `It seems to be
just block size~$2$ that is a problem.'  Perhaps those of us who usually deal
with larger blocks had simply not thought that it was worth while to 
investigate block size $2$ before the introduction of microarrays.

However, as we sketch in the next section, the problem is not block size~$2$ 
but very low average replication.  The proofs there are similar to those in 
this section;  they are given in more detail in \cite{rabalia, alia}.  
Once again, it turns out that these results are not all new.
The D-optimal designs for $v/(k-1)$ blocks of size~$k$ were given by 
Balasubramanian and Dey in \cite{baldey} in 1996---but their proof
uses a version of Theorem~\ref{thm:gaff} with the wrong value of the
constant.  The A-optimal 
designs for $v/(k-1)$ blocks of size~$k$ were published by Krafft and 
Schaefer in \cite{krafft} in 1997---but those authors are not blameless either,
because they apparently had not read \cite{tj1}!

Our best explanation is that agricultural statisticians are so familiar with
average replication being at least $3$ that when we saw these papers we decided
that they had no applicability and so forgot them.

\section{Very low average replication}

In this section we once again consider general block size~$k$.  A block 
design is connected if and only if its Levi graph is connected.  The Levi 
graph has $v+b$ vertices and $bk$ edges, so connectivity implies that 
$bk \geq b+v-1$; that is, $b(k-1) \geq v-1$.

\subsection{Least replication}

If $b(k-1) = v-1$ and the design is connected, then the Levi graph 
$\tilde G$ is a tree and the concurrence graph $G$ looks like those in
Fig.~\ref{fig:cgqueen}.  Hypergraph-theorists do not seem to have an
agreed name for such designs. 

For both D- and A-optimality, it turns out to be convenient to use the
Levi graph.   Since all the Levi graphs are trees,
Theorem~\ref{thm:gaff} shows that the D-criterion does not distinguish
among connected designs.


By 
Theorem~~\ref{thm:tjur}, $V_{ij} = \tilde R_{ij} \sigma^2$.  When
$\tilde G$
is a tree, $\tilde R_{ij}=2$ when $i$ and $j$ are in the same block; 
otherwise, $\tilde R_{ij} = 4$ if any block containing~$i$ has a
treatment in common with any block containing~$j$; 
and otherwise, $\tilde R_{ij} \geq 6$.  The queen-bee designs are the only
ones for which $\tilde R_{ij}\leq 4$ for all $i$ 
and $j$, and so they are the A-optimal designs.

The non-trivial eigenvalues of a queen-bee design are $1$, $k$ and $v$, with 
multiplicities $b-1$, $b(k-2)$ and $1$, respectively.  If the design is not a 
queen-bee design, then there is a treatment~$i$ that is in more than one block 
but not in all blocks. Thus vertex~$i$ forms a cutset for the concurrence graph
$G$ which is not joined to every other vertex of $G$.  Cutset Lemma~2 shows 
that $\theta_1<1$.  Hence the E-optimal designs are also the queen-bee designs.

\subsection{One fewer treatment}

If $b(k-1)=v$, then the Levi graph $\tilde G$ has $bk$ edges and $bk$ vertices,
and so it contains a single cycle, which must be of some even
length~$2s$. If $2\leq s \leq b$, then the design is binary; 
if $s=1$, then there is a single
non-binary block, whose defect is~$1$.  In this case, $k\geq 3$, because each
block must have more than one treatment.

For $2\leq s \leq b$, let $\mathcal{C}(b,k,s)$ be the class of designs
constructed as follows.  Start with a loop design for $s$ treatments. 
Insert $k-2$ extra treatments into each block.  The remaining $b-s$
blocks all contain the same treatment
from the loop design, together with $k-1$ extra treatments. 
Figs.~\ref{fig:star} and~\ref{fig:levistar} show the concurrence graph and
Levi graph, respectively, of a design in $\mathcal{C}(6,3,6)$.

For $k\geq 4$, the designs in $\mathcal{C}(b,k,1)$ have one treatment
which occurs twice in one block and once in all other blocks, with the
remaining treatments all replicated once.  The class $\mathcal{C}(b,3,1)$
contains all such designs, and also those in which the treatment in
every block is the one which occurs only once in the non-binary block.

\begin{thm}
If $b(k-1)=v$, then the D-optimal designs are those in $\mathcal{C}(b,k,b)$.
\end{thm}

\begin{pf}

%
%

The Levi graph $\tilde G$ is unicyclic, so its number of spanning
trees is maximized when the cycle has maximal length.
Theorem~\ref{thm:gaff} shows that the D-optimal designs are precisely
those with $s=b$.
\halmos
\end{pf}

\begin{thm}
\label{thm:abig}
If $b(k-1)=v$ then the A-optimal designs are those in $\mathcal{C}(b,k,s)$, 
where the value of $s$ is given in Table~\ref{tab:sval}.  
\end{thm}

\begin{table}[htbp]
\[
\begin{array}{cc|ccccccccccccc}
k & b & 2 & 3 & 4 & 5 & 6 & 7 & 8 & 9 & 10& 11 & 12 & 13\\
\hline
2 & & 2 & 3 & 4 & 5 & 6 & 7 & 8 & 4 & 4 & 4 & 3 \mbox{ or } 4& 3\\
3 & & 2 & 3 & 4 & 5 & 6 & 3 & 3 & 3 & 3 & 3 & 2 & 2\\
4 & & 2 & 3 & 4 & 5 & 3 & 2 & 2 & 2 & 2 & 2 & 2 & 2 \\
5 & & 2 & 3 & 4 & 5 & 2 & 2 & 2 & 2 & 2 & 2 & 2 & 2\\
6 & & 2 & 3 & 4 & 2 & 2 & 2 & 2 & 2 & 2 & 2 & 2 & 2
\end{array}
\]
\caption{Value of $s$ for A-optimal designs for $b(k-1)$ treatments in
$b$ blocks of size $k$: see Theorem~\ref{thm:abig}}
\label{tab:sval}
\end{table}

\begin{pf}
The Levi graph $\tilde G$ has one cycle, whose length is
$2s$, where $1\leq s\leq b$.  A similar argument to the one used at the start
of the proof of Theorem~\ref{thm:aopt} shows that this cannot be A-optimal
unless the design is in $\mathcal{C}(b,k,s)$.
If $s\geq 2$ or $k\geq 4$, then each block-vertex in the cycle has 
$k-2$ treatment-vertices attached as leaves; all other block-vertices
are joined to the same  single treatment-vertex in the cycle, and each has 
$k-1$ treatment vertices attached as leaves.  In $\mathcal{C}(b,3,1)$
the first type of design has a Levi graph like this, and the other
type has the same multiset of effective resistances between
treatment-vertices, because their concurrence graphs are identical.  
The following calculations use the first type.

Let $\mathcal{V}_1$ be the set of treatment-vertices in the cycle,
$\mathcal{V}_2$ the set of other treatment-vertices joined to blocks
in the cycle, and $\mathcal{V}_3$ the set of remaining
treatment-vertices. For $1\leq i\leq j\leq 3$, denote by
$\mathcal{R}_{ij}$ the sum of the pairwise resistances between
vertices in $\mathcal{V}_i$ and $\mathcal{V}_j$.

Put
\[
R_1 = \sum_{d=1}^{s-1}\frac{2d(2s-2d)}{2s} = \frac{s^2-1}{3}
\]
and
\[
R_2 = \sum_{d=0}^{s-1}\frac{(2d+1)(2s-2d-1)}{2s} = \frac{2s^2+1}{6}.
\]
Then $\mathcal{R}_{11} = sR_1/2$, $\mathcal{R}_{12} = s(k-2)(R_2+s)$, 
$\mathcal{R}_{13} = (b-s)(k-1)(R_1+2s)$,
$\mathcal{R}_{22} = s(k-2)(k-3) + s(k-2)^2[R_1 + 2(s-1)]/2$,
$\mathcal{R}_{23} = (b-s)(k-1)(k-2)(R_2+3s)$, and
$\mathcal{R}_{33} = (b-s)(k-1)(k-2) + 2(b-s)(b-s-1)(k-1)^2$.
Hence the sum of the pairwise effective resistances between
treatment-vertices in the Levi graph is $g(s)/6$, where
\[
g(s) = -(k-1)^2s^3 +2b(k-1)^2 s^2 -[6bk(k-1) -4k^2 +2k -1]s +c
\]
and $c=b(k-1)[12b(k-1) -5k-4]$.

If $s=1$ then the design is non-binary.  However, \[g(1)-g(2) =
(3k-9+6b)(k-1) -3,\] which is positive, because $k\geq 2$ and $b\geq
2$. Therefore the non-binary designs are never A-optimal.

Direct calculation shows that $g(2) > g(3)$ when $b=3$, and that
$g(2)> g(3)> g(4)$ when $b=4$.  These inequalities hold for all values
of $k$, even though $g$~is not decreasing on the interval $[2,4]$ for
large $k$ when $b=4$.

If $b=5$ and $k\geq 6$, then $g(3)>g(2)$ and $g(5)>g(2)$.  Thus the
local minimum of $g$ occurs in the interval $(1,3)$ and is the overall
minimum of $g$ on the interval $[1,5]$.

Differentiation gives
\[
g'(b) = b(k-1)[(b-6)(k-1) - 6] +4k^2 - 2k +1.
\]
If $g'(b)>0$ then $g$ has a local minimum in the interval $(1,b)$.
If, in addition, $g(3)>g(2)$, then the minimal value for integer $s$
occurs at $s=2$.  These conditions are both satisfied if $k=3$ and
$b\geq 12$, $k= 4$ and $b\geq 8$, $k\geq 5$ and $b\geq 7$, or $k\geq 9$
and $b\geq 6$.

Given Theorem~\ref{thm:aopt}, there remain only a finite number of
pairs $(b,k)$ to be checked individually to find the smallest value of
$g(s)$. The results are in Table~\ref{tab:sval}.
\halmos
\end{pf}

\begin{thm}
If $b(k-1)=v$, $b\geq 3$ and $k\geq 3$, then the E-optimal designs
are those in $\mathcal{C}(b,k,b)$ if $b\leq4$, and those in 
$\mathcal{C}(b,k,2)$ and $\mathcal{C}(b,k,1)$ if $b\geq 5$.
\end{thm}

\begin{pf}
If $2<s<b$ then the concurrence graph $G$ has a vertex which forms a 
vertex-cutset and which is not joined to all other vertices; moreover,  
$G$ has no multiple edges. Thus Cutset Lemma~2 shows that $\theta_1<1$.

Direct calculation shows that $\theta_1=1$ if $s=1$ or $s=2$.  For
$k\geq 4$, all contrasts between singly replicated treatments in the
same block are eigenvectors of the Laplacian matrix~$L$ with 
eigenvalue~$k$.  When $k\geq 3$ and $s=b$ the contrast
between singly and doubly replicated treatments has eigenvalue
$2(k-1)$. For $s=b$, a straightforward calculation shows that the
remaining eigenvalues of $L$ are
\[
k - \cos \left ( \frac{2\pi n}{b} \right) \pm 
\sqrt{(k-1)^2 - \sin^2 \left( \frac{2\pi n}{b}\right)}
 \]
for $1\leq n\leq b-1$.  The smallest of these is
\[k - \cos(2\pi/b) - \sqrt{(k-1^2) - \sin^2(2\pi/b)}:\]
this is greater than $1$ if $b=3$ or $b=4$, but less than $1$ if
$k\geq 3$ and $b\geq 5$.
\halmos
\end{pf}

\section{Further reading}

The Laplacian matrix of a graph, and its eigenvalues, are widely used, 
especially in connection with network properties such as connectivity,
expansion, and random walks. A good introduction to this material can
be found in the textbook by Bollob\'as~\cite{bollobas}, especially Chapters II
(electrical networks) and IX (random walks). Connection between the smallest
non-zero eigenvalue and connectivity is described in surveys by 
de Abreu~\cite{survey} and by Mohar~\cite{mohar}.
In this terminology, a version of Theorem~\ref{thm:eopt} is in
\cite{algconn}.

The basic properties of electrical networks can be found in textbooks of
electrical engineering, for example Balabanian and Bickart~\cite{wag}. A
treatment connected to the multivariate Tutte polynomial appears in Sokal's
survey~\cite{sokal}. Bollob\'as describes several approaches to the theory,
including the fact (which we have not used) that the current flow minimises
the power consumed in the network, and explains the interactions between
electrical networks and random walks in the network. See also Deo~\cite{deo}.

The connection with optimal design theory was discussed in detail by the 
authors in their survey~\cite{bcc09}. Further reading on optimal design 
can be found in John and Williams~\cite{JAJERW:book}, Schwabe~\cite{Schwab},
or Shah and Sinha~\cite{ss}. For general principles of experimental design,
see Bailey~\cite{rabbook}.

\paragraph{Acknowledgement}
This chapter was written at the Isaac Newton Institute for
Mathematical Sciences, Cambridge, UK, during the 2011 programme on
Design and Analysis of Experiments.

\end{document}